\magnification=\magstep0
\font \normtext=cmr10 scaled \magstephalf
\font \normmath=cmmi10 scaled \magstephalf
\font \normmathsym=cmsy10 scaled \magstephalf
\font \normmathops=cmex10 scaled \magstephalf
\font \normitalic=cmti10 scaled \magstephalf
\font \normslanted=cmsl10 scaled \magstephalf
\font \normfett=cmbx10 scaled \magstephalf
\font \normstrich=msbm10 scaled \magstep0
\font \normCPcpgross=cmr10 scaled \magstep1
\font \normCPcpklein=cmr7 scaled \magstep2
\font \kleintext=cmr7 scaled \magstephalf
\font \kleinmath=cmmi7 scaled \magstephalf
\font \kleinmathsym=cmsy7 scaled \magstephalf
\font \kleinmathops=cmex10 scaled \magstep0
\font \kleinitalic=cmti10 scaled \magstep0
\font \kleinslanted=cmsl10 scaled \magstep0
\font \kleinfett=cmbx7 scaled \magstephalf
\font \kleinstrich=msbm7 scaled \magstephalf
\font \kleinCPcpgross=cmr7 scaled \magstep1
\font \kleinCPcpklein=cmr5 scaled \magstep2
\font \klkleintext=cmr5 scaled \magstephalf
\font \klkleinmath=cmmi5 scaled \magstephalf
\font \klkleinmathsym=cmsy5 scaled \magstephalf
\let \klkleinmathops=\kleinmathops

\font \klkleinfett=cmbx5 scaled \magstephalf
\font \klkleinstrich=msbm5 scaled \magstephalf
\font \kleinsttext=cmr5 scaled \magstep0
\font \kleinstmath=cmmi5 scaled \magstep0
\font \kleinstmathsym=cmsy5 scaled \magstep0
\let \kleinstmathops=\kleinmathops

\font \kleinstfett=cmbx5 scaled \magstep0
\font \kleinststrich=msbm5 scaled \magstep0
\font \GRnormtext=cmr10 scaled \magstep2
\font \GRnormmath=cmmi10 scaled \magstep2
\font \GRnormmathsym=cmsy10 scaled \magstep2
\font \GRnormmathops=cmex10 scaled \magstep2
\font \GRnormitalic=cmti10 scaled \magstep2
\font \GRnormslanted=cmsl10 scaled \magstep2
\font \GRnormfett=cmbx10 scaled \magstep2
\font \GRnormstrich=msbm10 scaled \magstep2
\let \GRossCPcpgross=\GRnormtext
\font \GRossCPcpklein=cmr7 scaled \magstep3
\let \GRkleintext=\normCPcpklein
\font \GRkleinmath=cmmi7 scaled \magstep2
\font \GRkleinmathsym=cmsy7 scaled \magstep2
\let \GRkleinmathops=\kleinmathops

\font \GRkleinfett=cmbx7 scaled \magstep2
\font \GRkleinstrich=msbm7 scaled \magstep2
\let \GRklkleintext=\kleinCPcpklein
\font \GRklkleinmath=cmmi5 scaled \magstep2
\font \GRklkleinmathsym=cmsy5 scaled \magstep2
\let \GRklkleinmathops=\kleinmathops

\font \GRklkleinfett=cmbx5 scaled \magstep2
\font \GRklkleinstrich=msbm5 scaled \magstep2
\font \GZGRnormtext=cmr10 scaled \magstep4
\font \GZGRnormmath=cmmi10 scaled \magstep4
\font \GZGRnormmathsym=cmsy10 scaled \magstep4
\font \GZGRnormmathops=cmex10 scaled \magstep4
\font \GZGRnormitalic=cmti10 scaled \magstep4
\font \GZGRnormslanted=cmsl10 scaled \magstep4
\font \GZGRnormfett=cmbx10 scaled \magstep4
\font \GZGRnormstrich=msbm10 scaled \magstep4
\let \GZGRossCPcpgross=\GZGRnormtext
\font \GZGRossCPcpklein=cmr10 scaled \magstep3
\font \GZGRkleintext=cmr7 scaled \magstep4
\font \GZGRkleinmath=cmmi7 scaled \magstep4
\font \GZGRkleinmathsym=cmsy7 scaled \magstep4
\let \GZGRkleinmathops=\GRnormmathops

\font \GZGRkleinfett=cmbx7 scaled \magstep4
\font \GZGRkleinstrich=msbm7 scaled \magstep4
\font \GZGRklkleintext=cmr5 scaled \magstep4
\font \GZGRklkleinmath=cmmi5 scaled \magstep4
\font \GZGRklkleinmathsym=cmsy5 scaled \magstep4
\let \GZGRklkleinmathops=\kleinmathops

\font \GZGRklkleinfett=cmbx5 scaled \magstep4
\font \GZGRklkleinstrich=msbm5 scaled \magstep4
\def\leer{}%
\let\Schrstr=/%
\def\Sect#1{\S}%
%
%

\newcount\Ol%
\Ol=1%
\def\druckintanf{1}%
\def\druckintend{9999}%
\def\aktTeil{{\orgnorm\CPcp\/DIES SCHEINT NUR EIN \/TESTFILE ZU SEIN}}%
\def\aktAbschnitt{{\orgnorm Dies scheint nur ein Testfile zu sein}}%
\def\aktPargr{{\orgnorm\CPcp 0}}%
\def\aktSetc{{\orgnorm\CPcp 0}}%
\def\aktFile{Datei???}%
\def\lastFile{Datei???}%
\def\Pseudodatei{Datei???}%
\def\File"#1"{\global\def\aktFile{#1}%
\ifx\lastFile\Pseudodatei\global\let\lastFile=\aktFile\else\fi}%
\def\lastTeil~/{}%
\def\lastAbschnitt{{\orgnorm Erste Seite}}%
\newfam\itfam\newfam\bffam\newfam\slfam
\def\norm{\def\rm{\fam0\normtext}%
\textfont0=\normtext \scriptfont0=\kleintext \scriptscriptfont0=\klkleintext
\textfont1=\normmath \scriptfont1=\kleinmath \scriptscriptfont1=\klkleinmath
\textfont2=\normmathsym \scriptfont2=\kleinmathsym \scriptscriptfont2=\klkleinmathsym
\textfont3=\normmathops \scriptfont3=\kleinmathops \scriptscriptfont3=\klkleinmathops
\textfont\itfam=\normitalic \def\it{\fam\itfam\normitalic}%
\textfont\slfam=\normslanted \def\sl{\fam\slfam\normslanted}%
\textfont\bffam=\normfett
\scriptfont\bffam=\kleinfett\scriptscriptfont\bffam=\klkleinfett
                                   \def\bf{\fam\bffam\normfett}%
\let\normBBbb=\normstrich \let\kleinBBbb=\kleinstrich \let\kleinstBBbb=\klkleinstrich
\let\normzusatz=\normspez \let\kleinzusatz=\kleinspez \let\kleinstzusatz=\klkleinspez
\let\CPcpgross=\normCPcpgross \let\CPcpklein=\normCPcpklein
\baselineskip=15pt plus 1pt%
\lineskip=1pt%
\lineskiplimit=-2pt%
\setbox\strutbox=\hbox{\vrule height9.5pt depth4pt width0pt}%
\rm\hem=\Ol em}%
\def\klein{\def\rm{\fam0\kleintext}%
\textfont0=\kleintext \scriptfont0=\klkleintext \scriptscriptfont0=\kleinsttext
\textfont1=\kleinmath \scriptfont1=\klkleinmath \scriptscriptfont1=\kleinstmath
\textfont2=\kleinmathsym \scriptfont2=\klkleinmathsym \scriptscriptfont2=\kleinstmathsym
\textfont3=\kleinmathops \scriptfont3=\klkleinmathops \scriptscriptfont3=\kleinstmathops
\textfont\itfam=\kleinitalic \def\it{\fam\itfam\kleinitalic}%
\textfont\slfam=\kleinslanted \def\sl{\fam\slfam\kleinslanted}%
\textfont\bffam=\kleinfett
\scriptfont\bffam=\klkleinfett \scriptscriptfont\bffam=\kleinstfett
                                         \def\bf{\fam\bffam\kleinfett}%
\let\normBBbb=\kleinstrich \let\kleinBBbb=\klkleinstrich \let\kleinstBBbb=\kleinststrich
\let\normzusatz=\kleinspez \let\kleinzusatz=\klkleinspez \let\kleinstzusatz=\kleinstspez
\let\CPcpgross=\kleinCPcpgross \let\CPcpklein=\kleinCPcpklein
\baselineskip=12pt plus 1pt%
\lineskip=1pt%
\lineskiplimit=-2pt%
\setbox\strutbox=\hbox{\vrule height8.5pt depth3.5pt width0pt}%
\rm\hem=\Ol em}%
\def\gross{\def\rm{\fam0\GRnormtext}%
\textfont0=\GRnormtext \scriptfont0=\GRkleintext \scriptscriptfont0=\GRklkleintext
\textfont1=\GRnormmath \scriptfont1=\GRkleinmath \scriptscriptfont1=\GRklkleinmath
\textfont2=\GRnormmathsym \scriptfont2=\GRkleinmathsym \scriptscriptfont2=\GRklkleinmathsym
\textfont3=\GRnormmathops \scriptfont3=\GRkleinmathops \scriptscriptfont3=\GRklkleinmathops
\textfont\itfam=\GRnormitalic \def\it{\fam\itfam\GRnormitalic}%
\textfont\slfam=\GRnormslanted \def\sl{\fam\slfam\GRnormslanted}%
\textfont\bffam=\GRnormfett
\scriptfont\bffam=\GRkleinfett\scriptscriptfont\bffam=\GRklkleinfett
                                            \def\bf{\fam\bffam\GRnormfett}%
\let\normBBbb=\GRnormstrich \let\kleinBBbb=\GRkleinstrich \let\kleinstBBbb=\GRklkleinstrich
\let\normzusatz=\GRnormspez \let\kleinzusatz=\GRkleinspez \let\kleinstzusatz=\GRklkleinspez
\let\CPcpgross=\GRossCPcpgross \let\CPcpklein=\GRossCPcpklein
\baselineskip=18pt plus 1pt%
\lineskip=1pt%
\lineskiplimit=-3pt%
\setbox\strutbox=\hbox{\vrule height13pt depth5pt width0pt}%
\rm\hem=\Ol em}%
\def\GROSS{\def\rm{\fam0\GZGRnormtext}%
\textfont0=\GZGRnormtext \scriptfont0=\GZGRkleintext \scriptscriptfont0=\GZGRklkleintext
\textfont1=\GZGRnormmath \scriptfont1=\GZGRkleinmath \scriptscriptfont1=\GZGRklkleinmath
\textfont2=\GZGRnormmathsym \scriptfont2=\GZGRkleinmathsym \scriptscriptfont2=\GZGRklkleinmathsym
\textfont3=\GZGRnormmathops \scriptfont3=\GZGRkleinmathops \scriptscriptfont3=\GZGRklkleinmathops
\textfont\itfam=\GZGRnormitalic \def\it{\fam\itfam\GZGRnormitalic}%
\textfont\slfam=\GZGRnormslanted \def\sl{\fam\slfam\GZGRnormslanted}%
\textfont\bffam=\GZGRnormfett
\scriptfont\bffam=\GZGRkleinfett\scriptscriptfont\bffam=\GZGRklkleinfett
                                        \def\bf{\fam\bffam\GZGRnormfett}%
\let\normBBbb=\GZGRnormstrich \let\kleinBBbb=\GZGRkleinstrich \let\kleinstBBbb=\GZGRklkleinstrich
\let\normzusatz=\GZGRnormspez \let\kleinzusatz=\GZGRkleinspez \let\kleinstzusatz=\GZGRklkleinspez
\let\CPcpgross=\GZGRossCPcpgross \let\CPcpklein=\GZGRossCPcpklein
\baselineskip=28pt plus 2pt%
\lineskip=1pt%
\lineskiplimit=-5pt%
\setbox\strutbox=\hbox{\vrule height20pt depth8pt width0pt}%
\rm\hem=\Ol em}%
\let\orgnorm=\norm
\let\bleibklein=\klein
\newskip\extraskip
\newskip\vpc
\vpc=1pc plus 1pt minus 1pt%
\newdimen\hpc
\hpc=1pc%
\newdimen\hem
\let\noi=\noindent
\def\absatz{\par\noi}%
\def\orgneuzl{\hfill\break}%
\def\neuzl{\quad\orgneuzl}%
\def \mal{\mskip 1.5mu\mathord{\cdot}\nobreak\mskip 1.5mu}%
\def \notantiparallel{\setbox0=\hbox{$\uparrow\!\!\downarrow$}\dimen0=\wd0%
\mskip 1mu\hbox to \dimen0{$\uparrow$\hss$/$\hss$\downarrow$}\mskip 1mu}%
\def \notparallel{\setbox0=\hbox{$\uparrow\!\!\uparrow$}\dimen0=\wd0%
\mskip 1mu\hbox to \dimen0{$\uparrow$\hss$/$\hss$\uparrow$}\mskip 1mu}%
\let \Wkl=\angle
\def \oWkl{\hbox{\setbox1=\hbox{$\bigcirc$}\dimen1=\wd1%
\hbox to 0pt{\box1\hss}\hfill$\Wkl$ \hfill}}%
\def \qed{
\setbox1=\hbox{$\sqcap \!\!\!\! \sqcup$}
\dimen1=\wd1\hbox to \dimen1{\hfill}
\hfill\hbox to 0pt{\hss$\sqcap \!\!\!\! \sqcup$}\break
\hbox{~~~}\par\vskip-\baselineskip\medskip}
\def \qedbox{
\setbox1=\hbox{$\sqcap \!\!\!\! \sqcup$}
\dimen1=\wd1\hbox to \dimen1{\hfill}
\hfill\hbox to 0pt{\hss$\sqcap \!\!\!\! \sqcup$}\break
\hbox{~~~}\vskip-\baselineskip}
\def\qedbox{\hbox{$\sqcap \!\!\!\! \sqcup$}}
\def \CPcp{\CPcpklein\def \/##1{{\CPcpgross##1}}}%
\def\Zier#1{\ifmmode{\cal#1}\else$\cal#1$\fi}%
\def\restr#1_#2{\setbox1=\hbox{$#1$}\setbox2=\hbox{$_{#2}$}\dimen1=\ht1%
\dimen2=\dp2\box1\vrule width 0.4pt height \dimen1 depth \dimen2\mskip1mu \box2}%
\def\klrestr#1_#2{\setbox1=\hbox{$\scriptstyle{#1}$}%
\setbox2=\hbox{$\scriptstyle{_{#2}}$}\dimen1=\ht1%
\dimen2=\dp2\box1\vrule width 0.4pt height \dimen1 depth \dimen2\mskip1mu \box2}%
\def\displayrestr#1_#2{\setbox1=\hbox{$#1$}\setbox2=\hbox{$_{#2}$}\dimen1=\ht1%
\dimen2=\dp2\box1\vrule width 0.4pt height \dimen1 depth \dimen2\mskip1mu
\lower0.5ex\box2}%
\def\innerolra{\dimen1=\wd0\dimen2=\ht0\dimen3=0.5\wd0\advance\dimen3by1pt%
\hbox to \dimen3{\hss$\overleftarrow{\vrule width 0pt depth 0pt height \dimen2%
\hskip \dimen3\vrule width 0pt depth 0pt height \dimen2}$\hss}\hskip-1pt%
\hbox to 0pt{\hss\copy0\hss}\hskip-1pt%
\hbox to \dimen3{\hss$\overrightarrow{\vrule width 0pt depth 0pt height \dimen2%
\hskip \dimen3\vrule width 0pt depth 0pt height \dimen2}$\hss}}%
\def\drauf#1#2{{\setbox1=\hbox{#1}\setbox2=\hbox{#2}\dimen1=\wd1\fff%
\ifdim\dimen1<\wd2\dimen1=\wd2\else\fi\hbox to \dimen1{\hfill#1\hfill}%
\hskip-\dimen1\hbox to\dimen1{\hfill#2\hfill}}}%
\def\innerrechtsolra{\dimen1=\wd0\dimen2=\ht0%
\hbox to 0pt{\copy0\hss}%
\hbox to \dimen1{\hss$\overrightarrow{\vrule width 0pt depth 0pt height \dimen2%
\hskip \dimen1\vrule width 0pt depth 0pt height \dimen2}$\hss}}%
\hsize=19.5truecm %
\advance \hoffset by -2truecm%
\vsize=24.5truecm%
\advance \voffset by -1.25truecm%
\newdimen\halbhsize
\halbhsize=0.5\hsize
\newdimen\realhsize
\realhsize=\hsize
\newdimen\finalhsize
\finalhsize=\hsize
\newdimen\Zeilenlaenge
\def\setZeilenlaenge{\Zeilenlaenge=\hsize\advance\Zeilenlaenge by -\leftskip
\advance\Zeilenlaenge by -\hangindent}%
\def\ruek#1{\advance\leftskip by #1pc\advance\realhsize by -#1pc\setZeilenlaenge
}%
\newcount\rueknummer
\rueknummer=-1%
\def\noruek{\egroup\rueknummer=0\bgroup\GROSS\hangindent=1\hem\hangafter=01%
{}\setZeilenlaenge}%
\def\ruekeinf{\egroup\rueknummer=1\bgroup\gross\dimen1=1\hpc%
\advance\dimen1 by 1\hem\hangafter=01\hangindent=\dimen1\hskip 1\hpc%
\setZeilenlaenge}%
\def\ruekdopp{\egroup\rueknummer=2\bgroup\norm\dimen1=2\hpc
\advance\dimen1 by 1\hem\hangafter=01\hangindent=\dimen1\hskip 2\hpc%
\setZeilenlaenge}%
\def\ruekdreif{\egroup\rueknummer=3\bgroup\klein\dimen1=3\hpc
\advance\dimen1 by 1\hem\hangafter=01\hangindent=\dimen1\hskip 3\hpc
\setZeilenlaenge}%
\let \vwskip=\vskip%
\def \wskip#1\vpc{\par\penalty9999\vskip#1\vpc}%
\def \0#1{\ifnum#1=0\fff
\ifcase\rueknummer\wskip2\vpc\or\wskip4\vpc    \or
\global\def\lastAbschnitt{~}\or \wskip4\vpc\global\def\lastAbschnitt{~}\or
\wskip4\vpc\global\def\lastAbschnitt{~}\or\vfill\eject\global\def\lastAbschnitt{~}%
\else\fi\noi\noruek
\else\let\vwskip=\wskip\csname#1\endcsname\fi}%
\def \1{\ifcase\rueknummer\wskip1\vpc\or\wskip2\vpc\or\vwskip2\vpc\or\vwskip2\vpc\or\vwskip2\vpc
\else\fi\noi\ruekeinf}%
\def \2{\ifcase\rueknummer\wskip2\vpc\or\wskip1\vpc\or\vwskip1\vpc\or\vwskip1\vpc\or\vwskip1\vpc
\else\fi\noi\ruekdopp}%
\def \3{\ifcase\rueknummer\wskip2\vpc\or\wskip1\vpc\or\vwskip0.5\vpc\or
\vwskip0.5\vpc\or\vwskip0.5\vpc\else\fi\noi\ruekdreif}%
\def \bultem{\absatz\ifcase\rueknummer\noruek\or\ruekeinf\or\ruekdopp
\or\ruekdreif\or\ruekvierf\else\5\fi\hangafter=0%
\advance\hangindent by 1em\advance\realhsize by -1em\hskip-\rueknummer \hpc%
\hbox to 0pt{\hss\hbox to 1em{\hfil$\triangleright$\hfil}}}%
\def \mitem#1{\absatz\ifcase\rueknummer\noruek\or\ruekeinf\or\ruekdopp
\or\ruekdreif\or\ruekvierf\else\5\fi\hangafter=0%
\advance\hangindent by 1em\hskip-\rueknummer \hpc%
\hskip-1em
\setbox0=\hbox{#1~}\dimen1=\wd0\advance\realhsize by -1em%
\ifdim\dimen1<1em\dimen1=1em\else\fi\hbox to \dimen1{\box0\hfill}}%
\def\notem{\quad\orgnotem}
\def \orgnotem{\absatz\ifcase\rueknummer\noruek\or\ruekeinf\or\ruekdopp
\or\ruekdreif\or\ruekvierf\else\fi
\hangafter=0\hskip-\rueknummer \hpc}%
\def\asmitem#1{\setbox0=\hbox{#1~}\dimen1=\wd0\hbox to 2em{%
\ifdim\dimen1<2em\box0\hfill\else\hss\box0\fi}}%
%
\def\Kr#1{\ifx~#1\def\Krr##1{\widetilde{\mathchoice{\hbox{\normBBbb ##1}}%
{\hbox{\normBBbb ##1}}{\hbox{\kleinBBbb ##1}}{\hbox{\kleinstBBbb ##1}}}}\else
\mathchoice{\hbox{\normBBbb #1}}{\hbox{\normBBbb #1}}{\hbox{\kleinBBbb#1}}%
{\hbox{\kleinstBBbb #1}}\def\Krr{}\fi\Krr}%
\def\Aut~{\mathord{\rm Aut}^\sim\!}%
\def\8{\ss}%
\def\KlaStr{\raise0.5ex\hbox{$\scriptscriptstyle
($}'\raise0.5ex\hbox{$\scriptscriptstyle)$}}%
\def\)^o{)\hbox{\setbox0=\hbox{$)$}\dimen0=0.5\wd0\hbox to \dimen0{%
\hss$\buildrel\circ\over{\phantom{)}}$}}}%
\newdimen\pseudo%
\def \fffff{\global\advance\pseudo by 0pt}%
\def \fff{\hbox to 0pt{\hss}}%
\newcount\Teilnummer\Teilnummer=0%
\newcount\Abschnittsnummer\global\Abschnittsnummer=0%
\newcount\Paragraphennummer\global\Paragraphennummer=0%
\newcount\Satzetceteranummer\global\Satzetceteranummer=0%
\newcount\Formelnummer\global\Formelnummer=0 %
\newcount\Ziffernummer\global\Ziffernummer=0%
\newcount\Ordinalnummer\global\Ordinalnummer=0%
\def\Titel#1``#2''{\ifnum\Teilnummer=0\fff\global\Teilnummer=#1\fff\else
\global\advance\Teilnummer by 1\fff\ifnum\Teilnummer=#1\fff\else\errmessage{%
Zaehlerwirrwarr}\fi\fi\dimen1=0.5\hsize\vbox{
\leftskip=0pt plus \dimen1\rightskip=0pt plus \dimen1\vskip0.5cm%
\CPcp\noi#2\break\hbox{~~~}\vskip -1ex}
\global\def\aktTeil{{\CPcpgross TEIL \uppercase\expandafter{\romannumeral #1}}}
\global\Abschnittsnummer=0\fff\global\Paragraphennummer=0\fff\global\Satzetceteranummer=0\fff
\global\Formelnummer=0\global\Ziffernummer=0\global\Ordinalnummer=0}%
\def\Teil#1``#2''{\ifnum\Teilnummer=0\fff\global\Teilnummer=#1\fff\else
\global\advance\Teilnummer by 1\fff\ifnum\Teilnummer=#1\fff\else\errmessage{%
Zaehlerwirrwarr}\fi\fi\dimen1=0.5\hsize\vbox{\hrule height 1pt\vskip 1ex%
\leftskip=0pt plus \dimen1\rightskip=0pt plus \dimen1%
\centerline{TEIL \uppercase\expandafter{\romannumeral #1}}%
\CPcp #2\break\hbox{~~~}\vskip -1ex\hrule height 1pt}%
\global\def\aktTeil{{\CPcpgross TEIL \uppercase\expandafter{\romannumeral #1}}}
\global\Abschnittsnummer=0\fff\global\Paragraphennummer=0\fff\global\Satzetceteranummer=0\fff
\global\Formelnummer=0\global\Ziffernummer=0\global\Ordinalnummer=0}%
\def\Abschnitt#1``#2''{\ifnum\Abschnittsnummer=0\fff\global\Abschnittsnummer=#1\fff\fff\else
\global\advance\Abschnittsnummer by 1\fff\ifnum\Abschnittsnummer=#1\fff\else%
\errmessage{Zaehlerwirrwarr}\fi\fi {\bf
\ifnum\Abschnittsnummer=0{}\else Abschnitt \fi\ifcase\Abschnittsnummer
Einleitung                  \or A\or B\or C\or D\or E\or F\or G\or H%
\or I\or J\or L\or M\or N\or P\or Q\or R\or S\or T\or U\or V\or
W\or X\or Y\or Z\else \errmessage{Zaehlerwirrwarr}\fi}:~~#2
\global\def\aktAbschnitt{{\CPcpklein
\ifnum #1=0{}           \else (PART     \fi\ifcase#1
\ \or A\or B\or C\or D\or E\or F\or G\or H%
\or I\or J\or L\or M\or N\or P\or Q\or R\or S\or T\or U\or V\or
W\or X\or Y \or Z\else \ \fi)      }}%
\global\Paragraphennummer=0\fff\global\Satzetceteranummer=0\fff
\global\Formelnummer=0\global\Ziffernummer=0\global\Ordinalnummer=0}%
\def\Paragr#1{\ifnum\Paragraphennummer=0\fff\global\Paragraphennummer=#1\fff\else
\global\advance\Paragraphennummer by 1\fff\ifnum\Paragraphennummer=#1\fff\else
\ifodd\Paragraphennummer\fff\else\global\advance\Paragraphennummer by 1\fff\fi
\ifnum\Paragraphennummer=#1\fff\else
\errmessage{Zaehlerwirrwarr}\fi\fi\fi \Sect:\number\Paragraphennummer\fff
\global\def\aktPargr{{\CPcpklein #1}}%
\global\def\aktSetc{{\CPcpklein 0}}%
\global\Satzetceteranummer=0\fff\global\Formelnummer=0\global\Ziffernummer=0%
\global\Ordinalnummer=0}%
\def\Setc#1{\ifnum\Satzetceteranummer=0\fff\global\Satzetceteranummer=#1\fff\else
\global\advance\Satzetceteranummer by 1\fff\ifnum\Satzetceteranummer=#1\fff\else
\ifodd\Satzetceteranummer\fff\else\global\advance\Satzetceteranummer by 1\fff\fi
\ifnum\Satzetceteranummer=#1\fff\else
\errmessage{Zaehlerwirrwarr}\fi\fi\fi \number\Paragraphennummer.\number
\Satzetceteranummer\fff\global\Formelnummer=0\global\Ziffernummer=0%
\global\def\aktSetc{{\CPcpklein  #1}}%
\global\Ordinalnummer=0\ }%
\def\Formel#1{\if (#1#1\else
\ifnum\Formelnummer=0\fff\global\Formelnummer=#1\fff\else
\global\advance\Formelnummer by 1\fff\ifnum\Formelnummer=#1\fff\else
\errmessage{Zaehlerwirrwarr}\fi\fi
(\number\Formelnummer)\fi}%
\def\Ziffer#1{\ifnum\Ziffernummer=0\fff\global\Ziffernummer=#1\fff\else
\global\advance\Ziffernummer by 1\fff\ifnum\Ziffernummer=#1\fff\else
\errmessage{Zaehlerwirrwarr}\fi\fi
(\romannumeral\Ziffernummer)}%
\def\Nummer:#1.{\ifnum\Ordinalnummer=0\fff\global\Ordinalnummer=#1\fff\else
\global\advance\Ordinalnummer by 1\fff\ifnum\Ordinalnummer=#1\fff\else
\errmessage{Zaehlerwirrwarr}\fi\fi
(\number\Ordinalnummer.)}%
\def\Ordinal:#1.{\ifnum\Ordinalnummer=0\fff\global\Ordinalnummer=#1\fff\else
\global\advance\Ordinalnummer by 1\fff\ifnum\Ordinalnummer=#1\fff\else
\errmessage{Zaehlerwirrwarr}\fi\fi
\number\Ordinalnummer.}%
\def\Textskip #1.#2(#3.#4.#5.#6){\01\hrule height0pt\vtop{\bleibklein
\line{$\dotfill$}\vskip 6pt%
\centerline {\CPcp \/JUMP IN THE ENUMBERATION SYSTEM --- TEXT NOT NECCESSARILY
MISSING
}\vskip 6pt}\goodbreak\vbox{%
\global\Ziffernummer=0\global\Formelnummer=0\global\Ordinalnummer=0%
\global\Satzetceteranummer=#6\global\Paragraphennummer=#5%
\global\Abschnittsnummer=#4\global\Teilnummer=#3\fffff\bleibklein
\line{$\dotfill$}
\global\def\aktSetc{{\CPcpklein  #6}}%
\global\def\aktTeil{{\CPcpgross TEIL \uppercase\expandafter{\romannumeral #3}}}
\global\def\aktPargr{{\CPcpklein #5}}%
\global\def\aktAbschnitt{{\CPcpklein
\ifnum #4=0{}           \else   (PART   \fi\ifcase#4
\ \or A\or B\or C\or D\or E\or F\or G\or H%
\or I\or J\or L\or M\or N\or P\or Q\or R\or S\or T\or U\or V\or
W\or X\or Y\or Z\else \ \fi)      }}}}%
\def\dtTextskip #1.#2(#3.#4.#5.#6){\01\hrule height0pt\vtop{\bleibklein
\line{$\dotfill$}\vskip 6pt%
\centerline {\CPcp 
\/ES WURDE \/TEXT \"UBERSPRUNGEN UND DIE \/Z\"AHLER AUF
\Ref #1.#2.#5.#6(0) NEUEINGESTELLT
}\vskip 6pt}\goodbreak\vbox{%
\global\Ziffernummer=0\global\Formelnummer=0\global\Ordinalnummer=0%
\global\Satzetceteranummer=#6\global\Paragraphennummer=#5%
\global\Abschnittsnummer=#4\global\Teilnummer=#3\fffff\bleibklein
\line{$\dotfill$}
\global\def\aktTeil{{\CPcpgross TEIL \uppercase\expandafter{\romannumeral #3}}}
\global\def\aktSetc{{\CPcpklein  #6}}%
\global\def\aktPargr{{\CPcpklein #5}}%
\global\def\aktAbschnitt{{\CPcpklein
\ifnum #4=0{}           \else (PART     \fi\ifcase#4
\ \or A\or B\or C\or D\or E\or F\or G\or H
\or I\or J\or L\or M\or N\or P\or Q\or R\or S\or T\or U\or V\or
W\or X\or \or Z \else \ \fi)      }}}}%
\def\Ref #1.#2.#3.#4(#5){{%
\edef\roemTeilnummer{\expandafter\string\csname\romannumeral\Teilnummer\endcsname}%
\def\innBuchstzwo##1{\ifcase##1\#\or A\or B\or C\or D\or
E\or F\or G\or H\or G\or I\or K\or L\or M\or N\or O\or P\or Q\or S\or T\or
U\or V\or W\or X\or Y\or Z\else\errmessage{Zaehlerwirrwarr}\fi}%
\edef\Buchstzwo{\innBuchstzwo\Abschnittsnummer}\edef\?{?}%
\lowercase{\edef\Parameins{\expandafter\string\csname#1\endcsname}}%
\def\Nreins{#1}\count1=\Teilnummer\def\Pkteins{.}%
\def\Erstparam{#1}
\def\Paramzwei{#2}\def\Nrzwei{#2}\count2=\Abschnittsnummer\def\Pktzwei{.}%
\if\#\Paramzwei\def\Nrzwei{Einl}\else\fi
\def\Paramdrei{#3}\def\Nrdrei{#3}\count3=\Paragraphennummer\def\Pktdrei{.}%
\def\Paramvier{#4}\def\Nrvier{#4}\count4=\Satzetceteranummer
\def\Paramfunf{#5}\def\Nrfunf{(#5)}\count5=\Formelnummer\def\Nrnull{}%
\def\Figstring{Fig}%
\def\innerFigtest(##1.##2){##1}\edef\Figtest{\innerFigtest(#5.)}%
\ifx\Erstparam\?\vrule width 1.2ex height 1.1ex depth 0.1ex\else    
\ifx\Paramzwei\?\def\Nrnull{Teil }\def\Pkteins{}\def\Nrzwei{}\def\Pktzwei{}%
\def\Nrdrei{}\def\Pktdrei{}\def\Nrvier{}\def\Nrfunf{}\else
\ifx\Parameins\roemTeilnummer\def\Nreins{}\def\Pkteins{}\else\let\Buchstzwo=\ss
\count3=-2\fff\count4=-2\fff\count5=-2\fff\fi\ifx\Paramdrei\?%
\if\#\Paramzwei\def\Pktzwei{eitung}\else%
\def\Nrnull{}\def\Pktzwei{}\fi
\def\Nrdrei{}\def\Pktdrei{}\def\Nrvier{}\def\Nrfunf{}\else
\ifx\Paramzwei\Buchstzwo\def\Nrzwei{}\def\Pktzwei{}\else\count3=-2\fff\count4=2%
\fff\count5=-2\fff\def\Pktzwei{.}
\fi\ifx\Paramvier\?\def\Nrnull{\Sect.}\def\Pktdrei{}\def\Nrvier{}%
\def\Nrfunf{}\else\ifnum#3=\count3\def\Nrdrei{}\def\Pktdrei{}\else\count4=-2%
\fff\count5=-2\fff\fi
\ifx\Paramfunf\?\def\Nrfunf{}\def\Nrdrei{#3}\def\Pktdrei{.}%
\else\ifnum#4=\count4\def\Nrvier{}\else\count5=-2%
\def\Nrdrei{#3}\def\Pktdrei{.}\fi\fi\fi\fi\fi
\ifx\Nrvier\leer\ifx\Figtest\Figstring\def\Nrfunf{#5}\else\fi\else\fi
\Nrnull\Nreins\Pkteins\Nrzwei\Pktzwei\Nrdrei\Pktdrei\Nrvier\Nrfunf\fi}}%
\def\Arginner{\setbox2=\hbox{$2$}\dimen0=\wd0\dimen1=0.5\ht2\dimen2=\dimen1%
\advance\dimen1 by 0.2pt \advance\dimen2 by -0.2pt \dimen5=0.5\ht0\dimen6=\dp0%
\hbox to \dimen0{\vrule width 0.4pt height \dimen5 depth \dimen6 \hfill
\hbox to 0pt{\hss\box0\hss}\hbox to 0pt{\hss\box1\hss}\hbox to 0pt{\hss
\vrule width \dimen0 height \dimen1 depth -\dimen2\hss}\hfill
\vrule width 0.4pt height \dimen5 depth \dimen6}}%
\newbox\Druckam
\newcount\Std
\Std=\time
\divide\Std by 60%
\count1=\Std
\multiply\count1 by -60%
\advance \time by \count1%
\setbox\Druckam=\hbox to 0pt{\kleinsttext Druck am \number\day . \number\month . %
   \number\year\ um \number\Std :\number\time\ Uhr  \hss}%
\footline{\copy\Druckam \hss \normitalic\folio \hss%
\hbox to 0pt{\kleinsttext\hss\ifx\aktFile\lastFile\aktFile\else%
\lastFile/\aktFile\fi}}%
%
%
\def\morelinedisplaycalculations{\count1=\rueknummer\fff\ifnum\rueknummer<0\fff
\count1=0\fff\else\fi\dimen1=\realhsize\advance\dimen1by\leftskip
\dimen2=\count1\hpc\advance\dimen1 by-\dimen2\fff\ifnum\rueknummer<0\fff\else
\advance\dimen1 by -1\hem\fi\dimen2=\dimen1\advance\dimen2 by -2em%
\setbox0=\hbox{$(4)$}\dimen3=\wd0\advance\dimen3 by \leftskip
\advance\dimen2 by -\dimen3%
\lineskiplimit=3pt\advance\lineskip by 3pt}%
\def \twolinedisplay$$#1&#2&#3$${{\morelinedisplaycalculations
$$\vcenter{\hsize=\dimen1\hbox to \dimen1{$
\vcenter{\hsize=\dimen3\noi\hfill$#3$\hfill}\quad
\vcenter{\hsize=\dimen2\hbox to \dimen2{$\displaystyle{#1}$\hfill}%
    \hbox to \dimen2{\hfill$\displaystyle{#2}$}}\quad$}}$$}\ignorespaces}%
\def \threelinedisplay$$#1&#2&#3&#4$${{\morelinedisplaycalculations
$$\vcenter{\hsize=\dimen1\hbox to\dimen1{$
\vcenter{\hsize=\dimen3\noi\hfill$#4$\hfill}\quad
\vcenter{\hsize=\dimen2\hbox to \dimen2{$\displaystyle{#1}$\hfill}%
\hbox to \dimen2{\hfill$\displaystyle{#2}$\hfill}\hbox to \dimen2{\hfill
$\displaystyle{#3}$}}\quad$}}$$}\ignorespaces}%
\def \onelinedisplay$$#1&#2$${{\def\/{\egroup$\hfil$\displaystyle\bgroup}%
\morelinedisplaycalculations
$$\hbox to \dimen1{\hbox to \dimen3{\hskip\leftskip\hfil$#2$\hfil}\quad\hfil$
\displaystyle{#1}$\hfil\quad}$$}\ignorespaces}%
\def\oneormorelinedisplay$$#1&#2$${{\def\/{\egroup$\hfil$\displaystyle\bgroup}
\def\_/{\egroup$\egroup\hbox to\dimen2\bgroup$\displaystyle\bgroup}%
\morelinedisplaycalculations$$\vcenter{\hsize=\dimen1\hbox to \dimen1{$
\vcenter{\hsize=\dimen3\noi\hfil$#2$\hfil}\quad
\vcenter{\hsize=\dimen2\hbox to \dimen2{$
\displaystyle{#1}$}}\quad$}}$$}\ignorespaces}%
\def\Zeile#1{\neuzl\hbox to \realhsize{\let\/=\hfil#1}}%
%
\let\orgplus=\+%
\def\+{}%
\def\sw#1#2#3{\if #3a\swPlusKommaAnf\else\if #3e\swPlusKommaEnd
\else\errmessage{"sw,+"-Error}\fi\fi\liesvor}%
\def\swPlusKommaAnf{\bgroup\def\kOmma{,}\def\plUs{+}\catcode`\,=13
\catcode`\+=13
\swPlusKommaInner}%
\catcode`\,=13
\catcode`\+=13
\def\swPlusKommaInner{\def,{\mathbin{\kOmma}}\def+{\mathord{\plUs}}%
\def\liesvor##1##2{\if n##1\else\DruckErrmessage\fi
\if f##2\else\DruckErrmessage\fi
}}%
\def\swPlusKommaEnd{\egroup\def\liesvor##1##2{\if n##1\else\DruckErrmessage\fi
\if d##2\else\DruckErrmessage\fi}}%
\catcode`\,=12
\catcode`\+=12
\def\DruckErrmessage{\errmessage{"sw,+"-Error}}%
\let\+=\orgplus
%
\let\orgmakeheadline=\makeheadline
\advance \voffset by 2.5pc\advance\vsize by -1pc%
\def \makeheadline{\hbox to \finalhsize{\vbox to 0pt{\vss
\rightskip=0pt plus 0.5\hsize\fff
\leftskip=0pt plus 0.5\hsize\orgnorm\noi\hsize=\finalhsize
\centerline{\lastTeil~/~\lastAbschnitt}\vskip 22.5pt}}}%
%
%
\long\def \by{\par\egroup\hbox{~~~}}%
\norm
\def \vruhle{\vrule width 0.5pt}%
\def \inner{\dimen9=\dimen8\advance\dimen9 by -0.5pt}%
\def \Tabkast#1#2{\dimen8=#1\inner
\def\/##1${$\neuzl\line{\hfill$##1$}}%
\vtop{\vskip1pt\raggedright\bleibklein
\hsize=\dimen9\noi\setbox0=\hbox{#2}\ifdim\wd0<\dimen9%
\hbox to \dimen9{\hfill#2\hfill}\else#2\fi\vskip1pt}\vruhle}%
\newbox\Graphikbox
\setbox\Graphikbox=\vbox to 0pt{\vss\hrule depth 0pt height 0pt width 1truecm%
\vss}%
\newbox\GraphikPSbox
\setbox\GraphikPSbox=\vbox to 0pt{\vss\hrule depth 0pt height 0pt width 1truecm%
\vss}%
\newbox\Beschriftbox
\setbox\Beschriftbox=\vbox to 0pt{\vss\hrule depth 0pt height 0pt width 1truecm%
\vss}%
\def\bildrahmen#1x#2cm#3{\par \noi \offinterlineskip%
\dimen1=#1truecm\advance\dimen1 by -1truecm%
\line{\hfill\vrule\hbox{\vbox{\hrule width #1truecm\vskip #2truecm%
\vskip 1truecm\hbox{\hskip 0.5truecm\vbox{\hsize=\dimen1\par\noi%
{\bleibklein#3}}\hskip 0.5truecm}\vskip 2pt\hrule}}\vrule\hfill}%
\medskip}%
\def\orgbildbox#1x#2cm#3{\par \noi \offinterlineskip%
\dimen1=#1truecm\advance\dimen1 by -1truecm%
\line{\hfill\vrule\hbox{\vbox{%
\box\Graphikbox\box\Beschriftbox\hrule width #1truecm%
\vskip #2truecm\vskip 1truecm\hbox{\hskip 0.5truecm\vbox{\hsize=\dimen1\par\noi%
\bleibklein#3}\hskip 0.5truecm}\vskip 2pt\hrule width #1truecm}}\vrule\hfill}%
\medskip}%
\def\bildbox#1x#2cm#3{\par \noi \offinterlineskip%
\dimen1=#1truecm\advance\dimen1 by -1truecm%
\hbox to \realhsize{\hss\line{\hfill\vrule\hbox{\vbox{%
\box\Graphikbox\box\Beschriftbox\hrule width #1truecm%
\vskip #2truecm\vskip 1truecm\hbox{\hskip 0.5truecm\vbox{\hsize=\dimen1\par\noi%
\bleibklein#3}\hskip 0.5truecm}\vskip 2pt\hrule width #1truecm}}\vrule\hfill}%
}\medskip}%
\def\intextbox#1x#2cm#3{\hbox to \Zeilenlaenge{\hfill
\vrule\hbox{\vbox{\offinterlineskip\dimen1=#1truecm\advance\dimen1 by -1truecm%
\box\Graphikbox\box\Beschriftbox\hrule width #1truecm%
\vskip #2truecm\vskip 1truecm\hbox{\hskip 0.5truecm\vbox{\hsize=\dimen1\par\noi%
\bleibklein#3}\hskip 0.5truecm}\vskip 2pt\hrule width #1truecm}}\vrule\hfill}}
\def\pagebox#1{\par \noi \offinterlineskip%
\dimen1=\hsize\advance\dimen1 by -1truecm%
\line{\hfill\vrule\hbox{\vbox to \vsize{%
\box\Graphikbox\box\Beschriftbox\hrule width \hsize%
\vskip 18truecm\vfill\hbox{\hskip 0.5truecm\vbox{\hsize=\dimen1\par\noi%
\bleibklein#1}\hskip 0.5truecm}\vskip 2pt\hrule width \hsize}}\vrule\hfill}%
}
\def\BoxtextEnd.{\vrule width 0pt height 1pt depth 6pt.}%
\def\rechts{r}\def\oben{o}\def\mitte{m}%
\newdimen\xdimen
\newdimen\ydimen
\def\putbox nach (#1,#2) dieKoo (#3,#4) vonBox#5{%
\dimen1=#1truecm\dimen2=#2truecm\dimen3=\wd#5\dimen4=\ht#5%
\advance\dimen4 by \dp#5\def\Pardrei{#3}\def\Parvier{#4}%
\ifx\Pardrei\rechts \advance\dimen1 by -\dimen3\else\fi%
\ifx\Pardrei\mitte \advance\dimen1 by -0.5\dimen3\else\fi%
\ifx\Parvier\oben \advance\dimen2 by -\dimen4\else\fi%
\ifx\Parvier\mitte \advance\dimen2 by -0.5\dimen4\else\fi%
\vbox to 0pt{\hsize=\xdimen\vss\noi\hskip\dimen1\copy#5\vskip\dimen2}%
}
\hyphenation{struk-tur-er-hal-tend struk-tur-er-hal-ten-de
struk-tur-er-hal-ten-den mani-fold mani-folds
Struk-tur-in-va-ri-an-te Struk-tur-in-va-ri-an-ten
struk-tur-in-va-ri-an-te struk-tur-in-va-ti-an-ten}%
\def\Aehnlichkeitsto-{\"Ahn\-lich\-keits\-to\-}%
\let\Trennung=\-%
\def\Fktname#1({{\rm#1}(}%
\def\-{\ifmmode\Fktname\else\Trennung\fi}%
\newcount\druckanf%
\newcount\druckend%
\global\druckend=1%
\newcount\druckintnummer%
\global\druckintnummer=-1%
\newbox\muell%
\def\verschlucke{\global\setbox\muell=\hbox{\pagebody}\advancepageno
\setbox\muell=\hbox{.}}%
\output={\ifnum\pageno=\druckend\global\advance\druckintnummer by 1\fffff%
\global\druckanf=\ifcase\expandafter\druckintnummer\druckintanf\else0\fi\fffff%
\global\druckend=\ifcase\expandafter\druckintnummer\druckintend\else0\fi\fffff%
\global\advance\druckanf by-1\global\advance\druckend by1\fffff%
\else\fi%
\ifnum\pageno>\druckanf\ifnum\pageno<\druckend\plainoutput\else%
\verschlucke\fi\else\verschlucke\fi
\global\let\lastAbschnitt=\aktAbschnitt
\global\let\lastTeil=\aktTeil

\global\let\lastFile=\aktFile}%
\outer\long\def\ENDE{\message{ENDE erkannt}\par\egroup\hbox{~~~}%
\par\vfill\par\penalty-20000\ifnum\pageno<\druckend\else\deadcycles=0\fi%
\csname end\endcsname}%
\input epsf
\newbox\einzelpunkt
\setbox\einzelpunkt=\hbox{\normtext.}%
\newdimen\pkthoehe
\pkthoehe=\ht\einzelpunkt
\newdimen\pkthalbbreite
\newdimen\horzPSshift
\newdimen\downPSshift
\horzPSshift=0pt%
\downPSshift=0truecm%
\pkthalbbreite=0.5\wd\einzelpunkt
\newdimen\aktpos
\aktpos=0pt%
\newdimen\ydimen
\newdimen\xKOOanf
\newdimen\xKOOstep
\newdimen\xKOOend
\newdimen\yKOOanf
\newdimen\yKOOstep
\newdimen\yKOOend
\newdimen\xKOOsp
\newdimen\yKOOsp
\newdimen\xint
\newdimen\yint
\newcount\pktezahl
\newcount\pktenummer
\newcount\pkterest
\newcount\Graphiknr
\newdimen\GrosseKante
\newdimen\xKOOkorr
\catcode`\!=13
\catcode`\&=13
\catcode`\;=13
\catcode`\U=13
\def!#1{\rufz0#1}%
\def\rufz#1.#2#3#4#5#6{\yKOOanf=#1.#2#3#4#5\GrosseKante%
\xKOOanf=\xKOOend}%
\def\movexright#1{\advance#1by-\xKOOkorr\moveright#1}%
\def;#1{\if #1P\def\semik##1;{;}\else%
\if #1)\def\semik0){}\else
\def\semik##1.##2##3##4##5##6{%
\xKOOend=##1.##2##3##4##5##6\GrosseKante}\fi\fi\semik0#1}%
\defU{pt\global\pseudo=0pt\catcode`\#=6{}\uuhh}%
\def\uuhh#1PEN#2;{;}%
\def&#1{\kaufund0#1}%
\def\kaufund#1.#2#3#4#5#6{%
\yKOOend=#1.#2#3#4#5#6\GrosseKante%
\yKOOstep=\yKOOend\xKOOstep=\xKOOend%
\xint=\xKOOend\advance\xint by -\xKOOanf{}%
\yint=\yKOOend\advance\yint by -\yKOOanf{}%
\ifnum\xint<0{}\multiply\xint by -1{}\else\fi%
\ifnum\yint<0{}\multiply\yint by -1{}\else\fi%
\ifnum\xint>\yint{}\pktezahl=\xint\else\pktezahl=\yint\fi%
\xKOOsp=\xKOOend{}\yKOOsp=\yKOOend{}%
\ziehedielinie%
\xKOOanf=\xKOOsp{}\yKOOanf=\yKOOsp}%
\catcode`\!=12
\catcode`\&=4
\catcode`\;=12
\catcode`\U=11
\def\mitBeschriftung{\setbox\Beschriftbox=\vbox to 0pt\bgroup\offinterlineskip
\hrule height 0pt depth 0pt width \xdimen\vskip\ydimen}%
\def\soweitdieBeschriftung{\vss\egroup}%
\def\bildbox#1x#2cm#3{\par \noi \offinterlineskip%
\dimen1=#1truecm\advance\dimen1 by -1truecm%
\hbox to \realhsize{\hss\line{\hfill\vrule\hbox{\vbox{%
\box\Graphikbox\box\GraphikPSbox\box\Beschriftbox\hrule width #1truecm%
\vskip #2truecm\vskip 1truecm\hbox{\hskip 0.5truecm\vbox{\hsize=\dimen1\par\noi%
\bleibklein#3}\hskip 0.5truecm}\vskip 2pt\hrule width #1truecm}}\vrule\hfill}%
}\medskip}%

\def\initJozGraph(#1)#2x#3cm{%
\xdimen=#2truecm\ydimen=#3truecm%
\dimen7=\hsize%
\advance \dimen7 by -\xdimen%
\Graphiknr=#1%
\divide\dimen7 by2%
\advance\dimen7by\horzPSshift
\dimen8=-\xdimen\advance\dimen8 by\ydimen%
\advance\dimen8by\downPSshift%
\edef\JozGraphdatei{\ifcase\expandafter\Graphiknr\JozGraphdateien\else
keinname.pst\fi}%
\setbox\GraphikPSbox=\vbox to 0pt{\offinterlineskip%
\hbox{\hskip\horzPSshift\input\JozGraphdatei}\vss}%
\edef\altGraphikdatei{\ifcase\expandafter\Graphiknr\altGraphikdateien\else
pseudo.tex\fi}%
\setbox\Graphikbox=\vbox to 0pt\bgroup\offinterlineskip\vskip\ydimen%
\vskip0.5\pkthoehe%
\hrule width \xdimen height 0pt depth 0pt%
\xKOOkorr=0pt\GrosseKante=\xdimen\ifdim\xdimen<\ydimen\GrosseKante=\ydimen%
\xKOOkorr=\ydimen\advance\xKOOkorr by -\xdimen\divide\xKOOkorr by 2{}\else\fi%
\innerinitGraphik}%
\def\initGraphik(#1)#2x#3cm{%
\xdimen=#2truecm\ydimen=#3truecm%
\dimen7=\hsize%
\advance \dimen7 by -\xdimen%
\Graphiknr=#1%
\divide\dimen7 by2%
\advance\dimen7by\horzPSshift
\dimen8=-\xdimen\advance\dimen8 by\ydimen%
\advance\dimen8by\downPSshift%
\ifdim\downPSshift<0pt\edef\wche{-}\else\edef\wche{+}\fi%
\edef\Graphikdatei{\ifcase\expandafter\Graphiknr\Graphikdateien\else
keinname\fi}%
\setbox\GraphikPSbox=\vbox to 0pt{\offinterlineskip%
\epsfxsize=\xdimen\vskip\dimen8%
\hbox{\hskip\horzPSshift\epsfbox{\Graphikdatei}}\vss}%
\edef\altGraphikdatei{\ifcase\expandafter\Graphiknr\altGraphikdateien\else
pseudo.tex\fi}%
\setbox\Graphikbox=\vbox to 0pt\bgroup\offinterlineskip\vskip\ydimen%
\vskip0.5\pkthoehe%
\hrule width \xdimen height 0pt depth 0pt%
\xKOOkorr=0pt\GrosseKante=\xdimen\ifdim\xdimen<\ydimen\GrosseKante=\ydimen%
\xKOOkorr=\ydimen\advance\xKOOkorr by -\xdimen\divide\xKOOkorr by 2{}\else\fi%
\innerinitGraphik}%
\def\innerinitGraphik{%
\catcode`\!=13
\catcode`\&=13
\catcode`\;=13
\catcode`\U=13
\catcode`\#=14
\pseudo=%
\input\altGraphikdatei\vss\egroup}%
\def\initEpsGraphik(#1)#2x#3cm{%
\xdimen=#2truecm\ydimen=#3truecm%
\setbox\Graphikbox=\vbox{\hbox to \xdimen{#1}}}%
\def\ziehedielinie{%
\ifnum0=\xint{}\ifnum\yKOOanf<\yKOOend{}\linietypzwlf\else%
\ifnum0=\yint{}\linietypnull\else\linietypsex\fi\fi\else%
\ifnum0=\yint{}\ifnum\xKOOanf<\xKOOend{}\linietypdrei\else\linietypneun\fi\else%
\ifnum\xKOOanf<\xKOOend{}\ifnum\yKOOanf<\yKOOend{}\ifnum\xint>\yint{}%
\linietypzwei\else\linietypeins\fi\else\ifnum\xint>\yint{}%
\linietypvier\else\linietypfunf\fi\fi\else%
\ifnum\yKOOanf<\yKOOend{}\ifnum\xint>\yint{}\linietypzehn\else\linietypelf\fi%
\else\ifnum\xint>\yint\linietypacht\else\linietypsibn\fi\fi\fi\fi\fi}%
\def\schraeglinieFL{\gotonextsegmentFL}%
\def\schraeglinieST{\gotonextsegmentST}%
\def\gotonextsegmentFL{\ifdim\xKOOanf=\xKOOstep\else
\xint=\xKOOstep\advance\xint by -\xKOOanf%
\pktezahl=\xint\divide\pktezahl by 32768{}%
\ifnum\pktezahl=0{}\pktezahl=1{}\else\fi%
\nimmanfangssegment\performnextsegment\gotonextsegmentFL\fi}%
\def\gotonextsegmentST{\ifdim\yKOOanf=\yKOOstep\else
\yint=\yKOOstep\advance\yint by -\yKOOanf%
\pktezahl=\yint\divide\pktezahl by 32768{}%
\ifnum\pktezahl=0{}\pktezahl=1{}\else\fi%
\nimmanfangssegment\performnextsegment\gotonextsegmentST\fi}%
\def\nimmanfangssegment{\ifnum\pktezahl<30{}\else%
\advance\xKOOstep by \xKOOanf\divide\xKOOstep by 2%
\advance\yKOOstep by \yKOOanf\divide\yKOOstep by 2%
\divide\pktezahl by 2\advance\pktezahl by 1%
\nimmanfangssegment\fi}%
\def\performnextsegment{\dimen1=\xKOOanf\dimen2=\yKOOanf\dimen3=\xKOOstep%
\dimen4=\yKOOstep\pktenummer=-1{}\pkterest=\pktezahl\advance\pkterest by 1{}%
\let\nextpoint=\nextpointloop\nextpoint
\xKOOanf=\xKOOstep\xKOOstep=\xKOOend\yKOOanf=\yKOOstep\yKOOstep=\yKOOend}%
\def\nextpointloopend{}%
\def\nextpointloop{\advance\pktenummer by 1\advance\pkterest by -1{}%
\ifnum\pkterest<0{}\let\nextpoint=\nextpointloopend\else%
\dimen5=\dimen1\dimen6=\dimen2\dimen7=\dimen3\dimen8=\dimen4%
\multiply\dimen5by\pktenummer\multiply\dimen6by\pktenummer%
\multiply\dimen7by\pkterest\multiply\dimen8by\pkterest%
\advance\dimen5by\dimen7\advance\dimen6by\dimen8%
\divide\dimen5by\pktezahl\divide\dimen6by\pktezahl%
\dimen9=\dimen6%
\advance\dimen6by-\aktpos\advance\dimen6by\pkthoehe%
\advance\dimen5by-\pkthalbbreite%
\kern-\dimen6\movexright{\dimen5}\copy\einzelpunkt%
\aktpos=\dimen9{}%
\fi\nextpoint}%
\def\linietypnull{\dimen6=\yKOOanf\advance\dimen6by-\aktpos%
\advance\dimen6by\pkthoehe\dimen5=\xKOOanf\advance\dimen5by-\pkthalbbreite%
\kern-\dimen6\movexright{\dimen5}\copy\einzelpunkt\aktpos=\yKOOanf}%
\def\linietypeins{\schraeglinieST}%
\def\linietypzwei{\schraeglinieFL}%
\def\linietypdrei{\dimen6=\yKOOanf\advance\dimen6by-\aktpos%
\advance\dimen6by\pkthoehe\kern-\dimen6\movexright\xKOOanf\hbox{%
\vrule width \xint height \pkthoehe depth 0pt}\aktpos=\yKOOanf}%
\def\linietypvier{\schraeglinieFL}%
\def\linietypfunf{\invertyx\schraeglinieST}%
\def\linietypsex{\dimen6=\yKOOanf\advance\dimen6by-\aktpos\dimen5=\xKOOanf%
\advance\dimen5by-0.5\pkthoehe\kern-\dimen6\movexright{\dimen5}\hbox{%
\vrule width \pkthoehe height \yint depth 0pt}\aktpos=\yKOOend}%
\def\linietypsibn{\invertyx\schraeglinieST}%
\def\linietypacht{\invertyx\schraeglinieFL}%
\def\linietypneun{\dimen6=\yKOOend\advance\dimen6by-\aktpos%
\advance\dimen6by\pkthoehe\kern-\dimen6\movexright\xKOOend\hbox{%
\vrule width \xint height \pkthoehe depth 0pt}\aktpos=\yKOOend}%
\def\linietypzehn{\invertyx\schraeglinieFL}%
\def\linietypelf{\schraeglinieST}%
\def\linietypzwlf{\dimen6=\yKOOend\advance\dimen6by-\aktpos\dimen5=\xKOOanf%
\advance\dimen5by-0.5\pkthoehe\kern-\dimen6\movexright{\dimen5}\hbox{%
\vrule width \pkthoehe height \yint depth 0pt}\aktpos=\yKOOanf}%
\def\invertyx{\yKOOend=\yKOOanf\yKOOanf=\yKOOstep\yKOOstep=\yKOOend%
\xKOOend=\xKOOanf\xKOOanf=\xKOOstep\xKOOstep=\xKOOend}%
\def\prep(#1,#2)#3#4[#5#6]{\setbox8=\hbox{#5$#6$}\dimen8=0.#1\wd8\dimen9=\ht8%
\dimen7=\dp8\advance\dimen9 by \dimen7\dimen7=0.#2\dimen9%
\vbox to \dimen7{\if#4v\vss\else\fi\vbox to \dimen9{\offinterlineskip\vfil%
\hbox to \dimen8{\if#3>\hss\else\fi\box8\if#3<\hss\else\fi}\vfil}%
\if#4^\vss\else\fi}}%
\def\initaltGraphik(#1)#2x#3cm{%
\xdimen=#2truecm\ydimen=#3truecm%
\Graphiknr=#1%
\edef\altGraphikdatei{\ifcase\expandafter\Graphiknr\altGraphikdateien\else
keinname\fi}%
\setbox\Graphikbox=\vbox to 0pt\bgroup\offinterlineskip\vskip\ydimen%
\vskip0.5\pkthoehe%
\hrule width \xdimen height 0pt depth 0pt%
\xKOOkorr=0pt\GrosseKante=\xdimen\ifdim\xdimen<\ydimen\GrosseKante=\ydimen%
\xKOOkorr=\ydimen\advance\xKOOkorr by -\xdimen\divide\xKOOkorr by 2{}\else\fi%
\innerinitaltGraphik}%
\def\innerinitaltGraphik{%
\catcode`\!=13
\catcode`\&=13
\catcode`\;=13
\catcode`\U=13
\catcode`\#=14
\pseudo=%
\input\altGraphikdatei\vss\egroup}%
\def\|{\l}%
\def\7{\L}%
\let\altkomma=\,%
\def\,#1{\ifmmode\altkomma#1\else\c#1\fi}%
 
\hsize=12truecm\dimen1=\hsize\hsize=1.22\dimen1\fffff
\realhsize=\hsize\finalhsize=\hsize%
\vsize=19.5truecm\dimen1=\vsize\vsize=1.22\dimen1\fffff
\advance\vsize by -0.8truecm%
\def\4{\global\let\klein=\bleibklein\3\global\let\klein=\orgnorm\rueknummer=4{}}%
\def\ruekvierf{\global\let\klein=\bleibklein\ruekdreif\global\let\klein=\orgnorm\rueknummer=4{}}%
\def\Sect#1{}%
\let\Dach=\^%
\def\^#1{\ifx c#1\v#1\else
\ifx C#1\v#1\else
\ifx s#1\v#1\else
\ifx S#1\v#1\else
\ifx z#1\v#1\else
\ifx Z#1\v#1\else\Dach#1\fi\fi\fi\fi\fi\fi}
\let\klein=\norm
\def\norm{\klein\sl\def\sl{\it}}%
\Ol=0\fffff
\hpc=0pt%
\advance \hoffset by 2truecm%
\hfuzz 3.3pt%
\maxdeadcycles=66
\let\ttt=\tt\def\tt{\ttt\rightskip=0pt plus 0.88\hsize}%
\def\altGraphikdateien{pseudo.tex}%
\def\JozGraphdateien{pseudo.pst
\or figur_q.pst
\or figur_tt.pst
\or figur_z.pst
\or figur_w.pst
\or figur_u.pst
\or figur_r.pst
\or figur_ss.pst
\or figur_y.pst
\or figur_v.pst
\or keinname.pst}%
\def\Graphikdateien{%
\or 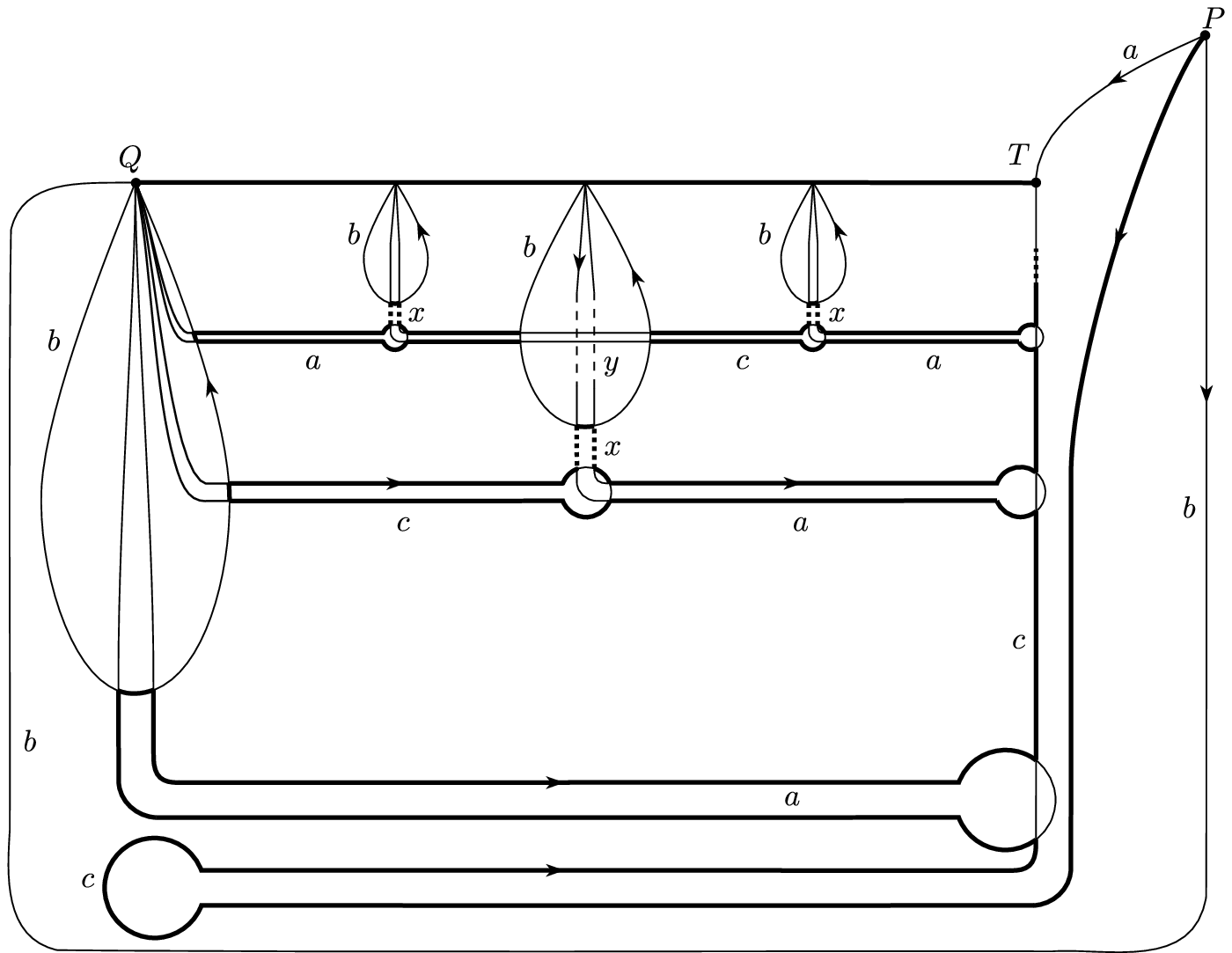%
\or T3undT_d.ps%
\or 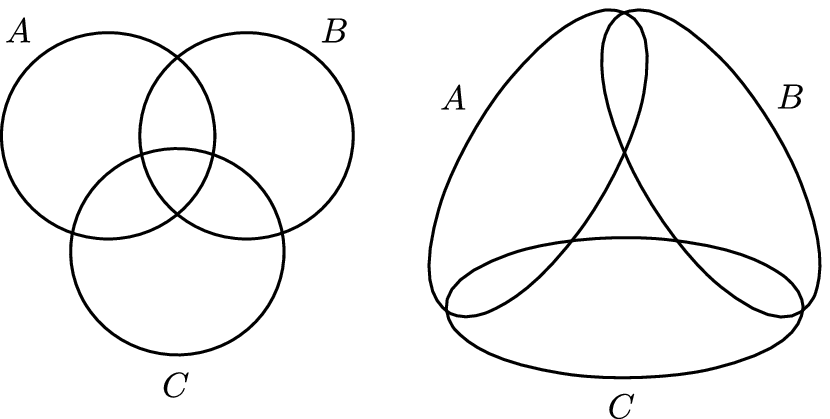%
\or revofall.ps%
\or keinname.tex}%
\rueknummer=0\fff
\let\makeheadline=\orgmakeheadline
\headline{\ifnum\pageno=1 \centerline{~~~}\else\centerline{
\CPcp \/D.\ \/REPOV\^S, W.\ \/ROSICKI, \/{\^Z}.\ \/VIRK \& \/A.\ \/ZASTROW}\fi}%
{\ruekdreif\norm\rm{
}\pageno=0\hbox{~}\vfill\eject{%
\vskip-3\baselineskip
\File"3c0l-preprint"%
\00\Titel1``\line{\klein Gda\'nsk \& Ljubljana $\dotfill$ Spring 2012}\break\line{\hfil}\break
 On Minc' sheltered middle path ''%
\vskip-1.5\baselineskip
\3\centerline{Du\^san Repov\^s, Witold Rosicki, \^Ziga Virk and Andreas Zastrow}\medskip
\4\hbox{~}\footnote{}{{\bleibklein\vskip-1em\noi
2010 Mathematics Subject Classification. Primary: 54D05, 54C05; Secondary: 54C50.
\neuzl Key words and phrases. Sheltered middle path, Helly type theorem, topologist's sine curve, winding number.}}

\medskip
\4\leftskip=1pc\rightskip=1pc%
{\bf Abstract: } 
This paper shows that a construction, which was
           introduced by Piotr Minc in connection with a problem
           that came from Helly-type theorems and that allows
           to replace three PL-arcs with a ``sheltered middle path", can
           in the case of general (non-PL) paths result in the
           topologist's sine curve.
\1{\bf\Paragr{1} . } Introduction:
Motivation via Helly-type-theorems and the three-colored sphere-problem
\3{\bf\Setc{1} Motivation}\neuzl
Helly's two-page theorem [Hlly1] from 1923 (cf.\ \Ref I.\#.1.3(?)), which was originally
a theorem from geometry about intersections of convex bodies, has given rise to
many generalizations and variations (see e.g.\ [DGK],
[Eck] \& [Wng], [Farb] \& [Flstd] for more recent papers
not cited in these survey articles). Among them are topological Helly theorems
(cf.\ \Ref I.\#.1.4(?)) which are
analogously structured implications where Helly's
original condition ``convex" is in the assumptions and the conclusions are replaced
by appropriate
topological conditions. In two dimensions Helly's theorem
implies, that if in a family of convex sets each three intersect,
then the entire family has a nonempty intersection. The first one
to prove topological Helly theorems was Helly himself in [Hlly2].
While his proof in general uses Vietoris homology, Helly in \S2
of [Hlly2] also
considers the two-dimensional case separately
(rather as a motivation, than as the initial step for some induction)
and works with simply connected sets.
A fairly natural picture that one has
of an intersection of three sets in the plane (cf.\ \Ref I.\#.1.1(Fig.a))
compared to an analogous picture for four sets
must have lead to the question, whether
 one can imagine three simply connected sets really
intersecting with connected intersections
in such a different way than pictured in \Ref I.\#.1.1(Fig.a).
In 1956 Moln\'ar asserted (quoted below as \Ref I.\#.1.5(?)) a statement
in this sense as a ``known theorem''. One
year later, in a paper [Mol2] that unfortunately only
appeared in Hungarian, he used this assertion as a startpoint
of an induction to prove a theorem that attempts to generalize
Helly's theorem for two dimensions so that he reduces the assumption
``simply connected" for the triple-intersections
to ``non-empty". He calls this theorem a ``two-dimensional Helly-theorem",
while later it got referred to as ``Moln\'ar's Theorem" (e.g.\ [Breen])
or Helly-Moln\'ar Theorem (e.g.\ [Bgt] in Rem.18
from \S4.3). It was stated in this form in that form
(without the additions ``path-connected")
in [DGK].\S4.11, and also quoted
in the Introduction of [Breen]. However,
in this form it was disproved by [KR\^Z].\S5 
together with Thm.1 of [Breen].

Apparently the first one to note that there was a problem
with these theorems was Bogatyi in [Bgt]. While his paper
is mainly written to close a gap in the proof of
Topological Helly Theorem for singular homology theory,
he also considers many related statements. Molnar's footnote
has also attracted his attention: Apparently he does
not believe that it is true as stated by Moln\'ar, but
he conjectures only (Hyp.2 in \S4.3 in [Bgt]) that the analogous version where
``connected" is replaced by ``path-connected" is true,
and he proves this hypothesis only in the special case
of Peano continua (it appears as Corl. 4.2 on page 397, where ``function"
   in this context has the very special meaning explained in the
   paragraph after the end of the proof of Prop.1 on page 367, but this statement is
best phrased in the pen-penultimate sentence
of the abstract).
Another weakened
version of Bogatyi's hypothesis, where one does
not need the assumption of local pathwise connectivity, and instead replaces
the condition ``connected" by ``pathwise connected" in the assumption, but
not in the conclusion, was proved
 by [KaRe], and this result was strengthened
independently in the papers by Tymchatyn and Valov ([TyVa].Prop.1.2) and
Minc ([Minc].Corl.1.2), by proving the full hypothesis of Bogatyi.
In his proof of
this version of Moln\'ar's footnote, Minc proposes the construction of a
``sheltered middle path''.\par
At the same time while Minc was working on
his paper we have for some time been building on the hypothesis stated as
Question \Ref I.\#.1.2(?) below:
\initGraphik(3)10x8cm
\topinsert\bildbox10x6cm{{\bf Figure a} (belonging to
\Ref I.\#.1.1(?)):               
The left figure is the classical Venn-diagram
  of three intersecting sets. The left and right figure together
  correspond to what one is tempted to believe to be the only other
  possibility for three simply connected sets in the plane
  to intersect with connected pairwise intersections, and this image explains
  Moln\'ar's assertion \Ref I.\#.1.5(?).
However cf.\ \Ref I.\#.1.6(?)\BoxtextEnd.}\endinsert\noi
Although [TyVa] and Minc' results solve the
problem of Moln\'ar's footnote, we are treating the three-colored
sphere-problem as a problem which is of interest on its own. The method of a
sheltered middle path solves the three-colored
sphere problem in the case when the three
bi-colored meridians are mapped into the plane without double points. It is not
difficult to see (and it is explained in the paper below, see
Lemma \Ref I.\#.2.4(?)),
that Minc' method extends to the case where the images of the three
bi-colored meridians have only finitely many intersections and
self-intersections. However, the goal of this paper is to demonstrate that
Minc' method of a sheltered middle path apparently cannot solve
the three-colored sphere problem in general.

\3{\bf  The main result.}
{\sl There exists a continuous map of the 3-colored sphere  into the plane such that the
sheltered middle-path does not give a topological path but only a topologist's sine curve.}
\3{\bf\Setc{2}  Question } {\sl (The three-colored sphere-problem):}\mitem{%
\it In Brief: }
      If three symmetric regions between the north and the south pole
      of a \hbox{2-sphere} $S^2$ are colored by different colors and this sphere
      is in an arbitrary continuous way mapped into $\Kr R^2$,
      must there always exist a path between the images
      of the north and of the south pole, each of whose
      points has preimages of all three colors?\notem
{\it More precisely stated,}\/ this is the following question:\neuzl
Consider a two-dimensional sphere $S^2$ and let $N$ and $S$ denote its north and
south pole, respectively. Connect $N$ and $S$ by three different meridians (paths) thus
subdividing $S^2$ into three different subsets $X_1$, $X_2$ and $X_3$ bounded by the
chosen meridians. We think of these subsets as closed, i.e., each of the three meridians
belongs to two of these subsets and  $X_1 \cap X_2 \cap X_3 =\{N,S\}$. Given an arbitrary
continuous map $f: S^2 \to \Kr R^2$, the points $f(N)$ and $f(S)$  lie in
$f(X_1)\cap f(X_2)\cap f(X_3)$. Is there always
a path between $f(N)$ and $f(S)$ which is
contained in $f(X_1)\cap f(X_2)\cap f(X_3)$?
\2{\bf\Setc{3} Theorem } {\it (Helly's original theorem { \rm  [Hlly1]):}\neuzl
%
Suppose that
$
    X_1,X_2,\ldots,X_n
$
is a finite collection of convex subsets of $\Kr{R}^d$, where $n > d$.
If the intersection of every $d + 1$ of these sets is nonempty, then the
entire collection has a nonempty intersection:
$
    \bigcap_{j=1}^n X_j\neq\emptyset.
$
\smallskip
\3{\bf\Setc{4} Remark: }
``Topological Helly Theorems" are statements, where the
convexity-assumptions for the sets $X_i$ are replaced by corresponding
topological assumptions (e.g.\ acyclicity with respect to some
homology theory), and they usually assert
similar assertions for the overall intersection.
\smallskip
\3{\bf\Setc{5} Moln\'ar's assertion from 1956 { \rm  ([Mol1]):}} (the parentheses ``(pathwise)" do not appear in Moln\'ar's text)
If $A$, $B$ and $C$ are three simply connected sets in the plane
and all three pairwise intersections are nonempty
and  (pathwise) connected,
and if the intersection of all three sets is nonempty, then
it is (pathwise) connected.
\vskip0.3\hsize
\3{\bf\Setc{6} Proven and disproven versions of Moln\'ar's assertion: }
\bigskip

{
\offinterlineskip
\def\vlong{\vrule height2.75ex depth1.25ex width 0.6pt}%
\def\vshrt{\vrule height1.75ex depth1.25ex width 0.6pt}%
\tabskip=0pt
\halign{
                                             #&#&\hfil#\hfil\tabskip=1em &#\vrule&
\hfil#\hfil &\vrule # &
\hfil #\hfil &#\hfil\vrule \tabskip=0pt \cr
\noalign{\hrule height 0.6pt}
\vlong&&Condition for intersections && Conclusion for&& Status of assertion&\cr
\vlong&&$A\cap B, A\cap C$ and $B\cap C$ &&   $A \cap B \cap C$ &&   &\cr
\noalign{\hrule}
\vlong(a) &&connected &&  simply connected && disproved by [KR\^Z] &\cr   
\vlong(b) &&path-connected && simply connected && conjectured by [Bgt] &\cr
\vlong(c)  &&connected (for locally \hfill &&locally path-connected\hfill&& proved by [Bgt]&\cr
\vshrt     && \hfill  path-connected $A, B\& C$)&& \hfill and simply connected&&&\cr
\vlong(d) &&path-connected && connected && proved by [KaRe] &\cr
\vlong(e) &&path-connected && path-connected && proved by [Minc]  &\cr
\vlong (f) &&path-connected&&simply connected&& proved by  [TyVa] &\cr
\noalign{\hrule height 0.6pt}
}}
\smallskip\noi Note that the discrepancy in the second column between
the statement of \Ref I.\#.1.5(?) and all lines apart from (e) can
easily be explained: If a subset of the plane
has non-trivial $\pi_1$, this is also unveiled  by a simple
closed non-nullhomotopic curve.
If such a curve can be found
in the intersection of three simply connected sets, the interior
of this curve must have belonged to each of the sets
and thus to the intersection, contradicting the assumption
that the curve was non-nullhomotopic. Thus the non-trivial claim
of the conclusion in lines (b)--(d) \& (f) is the triviality of $\pi_0$
(i.e.\ that the space is path-connected), not that of
$\pi_1$.
\3{\bf\Setc{7} How a solution of the three-colored sphere-problem
   implies Moln\'ar's assertion in the version (b) \& (d)--(f):}\neuzl
Let $P$ and $Q$ be any two points in $A\cap B\cap  C$.
Since $A\cap B$,  $B\cap C$ and $A\cap C$ are path-connected
we obtain three paths connecting $P$ with $Q$, one
in each of the intersections. The assumed simple connectivity
of $A$, $B$ and $C$  implies that each two of these
three paths are relatively homotopic. Hence we can define a map on a three-colored
sphere as described in Question \Ref I.\#.1.2(?) by the following rules:
\bultem the three paths from the intersection of two of the three
   sets give the maps  of the three bi-colored meridians.
\bultem the three homotopies are used to define the maps on the three
   one-colored regions.\notem
Since the definitions on the boundaries of one-colored regions coincide  we obtain a
continuous map of the entire sphere. The definition is consistent
with the coloring in the sense that each color is entirely mapped
into one of the sets $A$ or $B$ or $C$. Thus the existence of a
three-colored path between the image of the north and the image
of the south pole (which get mapped to $P$ and $Q$) implies
that this path lies in $A\cap B \cap C$. Consequently,  $P$ and $Q$
belong to the same path-component of $A\cap B\cap C$.

\3{\bf\Setc{8} Definition: } Let $f \colon S^1\to \Kr R^2$ denote a loop in the
plane and suppose $x\notin f(S^1).$ Given any point $y\in S^1$ define $\alpha_y$
as the positively oriented circular loop with center at $x$ and basepoint at $f(y)$.
The {\it``winding number''} of the point $x$ with respect to the loop $f$ is the
integer $k$ such that $[f]=[\alpha_y]^k\in \pi_1(\Kr R^2-\{x\},f(y))$. It turns out
that the definition is independent of the choice of point $y$.

Note that if $x$ has a nonzero winding number with respect to the loop $f$ then every
nullhomotopy of $f$ contains $x$ in its image.

\3{\bf\Setc{9} Remark: } The winding number is a well known tool which appears in
many areas of mathematics. In the case of a reasonably well behaved loop (i.e., a loop
allowing the construction below) it can be computed in the following way. Choose any
differentiable arc $p\colon [0,\infty)\to \Kr R^2$  satisfying the following conditions:
\mitem{(a)} $p(0)=x$, $|p'(t)|\neq 0 \quad \forall t$, and $p$ diverges to infinity,
i.e., $\lim_{t\to\infty}d(x,p(t))=\infty$;
\mitem{(b)} $p$ has finitely many intersections with loop $f$, none of which is a
self-intersection of $f$;
\mitem{(c)} every intersection of $p$ and $f$ is transversal, i.e., the loop $f$ is
locally differentiable at every intersection and the corresponding direction vectors
$p'$ and $f'$ are linearly independent at every intersection. The signature of the
intersection is defined to be $+1$ if the oriented angle between $p'$ and $f'$ is
between $0$ and $\pi$. If the oriented angle between $p'$ and $f'$ is between $\pi$
and $2\pi$  the signature of the intersection is $-1$.
\notem
If such an arc can be found then the winding number of $x$ with respect to $f$ can be
computed as the sum of the signatures of all intersections between the loop $f$ and the arc $p$.

\1{\bf\Paragr{2} . } The tame case: the construction of a three-colored path as  a sheltered middle path
                                                       in the PL-case}

\3An analysis of Minc' papers shows that in his
combinatorial proofs in Section 2  he was actually working
with a stronger version of sheltered than what he defined
on the first page of his paper. In the sequel we shall call the
definition from the beginning of his paper {\it``weakly sheltered",}\/
the definition that he was actually using {\it``strongly sheltered",}\/
and an intermediate version that would suffice to prove the
three-colored sphere-problem we shall call just {\it``sheltered"}.
Our counterexample from Sect.\Ref I.\#.3.?(?) will also be for this version
of shelteredness.

Given paths $a,b\colon [0,1]\to \Kr R^2$ the inverse path of $a$ will be denoted by
$a^-(t):=a(1-t)$; the concatenation of paths $a$ and $b$ that satisfy $a(1)=b(0)$
will be denoted by $a*b$.  If $a$ and $b$ are paths with $a(0)=b(0)$ and $a(1)=b(1)$, then
the loop determined by paths $a~\&~b$ is $a^- * b$ or $b^- * a$. Such loop is determined
up to an orientation.

\3{\bf\Setc{1} Definition: }
Let $a,b,c\colon [0,1]\to \Kr R^2$ be paths from $A$ to $B$, where $A\neq B$ are  points
in the plane. A point in the union of their traces is called:
\mitem{(a)} {\it``weakly sheltered'',}\/ if it lies either on at least two
of the three traces,
     or on  just one of the traces, but in the bounded
     component with respect to the union of the remaining two traces;
\mitem{(b)} {\it``sheltered'',}\/ if it lies either on at least two of three traces,
     or on  just one of the traces, but has non-zero
     winding number with respect to the loop determined by the remaining two paths;
\mitem{(c)} {\it``strongly sheltered'',}\/ if it lies either on at least two of three traces,
or on  just one of the traces, but has an odd
     winding number with respect to the loop determined by the remaining two paths.

The loop determined by two paths is unique up to orientation which does not affect the
properties of having even, odd or non-zero winding number.

\3{\bf\Setc{2} Lemma: }
Let $a,b$ and $c$ be three PL-paths in the plane from $A$ to $B$ which have only
general position intersections and self-intersections.
\mitem{(a)} If $P$ is an intersection of two different paths then exactly two of the four
incoming segments at $P$ are strongly sheltered. Furthermore, the two strongly sheltered
segments are not contained in the same path.
\mitem{(b)} If $P$ is a self-intersection of one of the paths then either none or all four
of the incoming segments at $P$ are strongly sheltered.
\mitem{(c)} If $P$ is either $A$ or $B$, then either one or all three of the three incoming
segments at $P$ are strongly sheltered.

\3{\it  Proof: } All three claims are proved using the winding number analysis of Remark
\Ref I.\#.1.9(?), hence the paths used in the proof satisfy appropriate conditions required
in that Remark. Note that the traces of the three PL-paths define a graph whose vertices
are induced by the PL-structure of paths and by their intersection. It is easy to see
that the strong shelteredness is an invariant of open edges (segments) of such a graph.
\mitem{(a)} Two points on different incoming segments of the same path can be connected
by a differentiable arc which transversely intersects the second path of the intersection
exactly once, hence changing the parity of the appropriate winding numbers of original
points. Consequently,  exactly one of the two incoming segments of each path is strongly
sheltered.
\mitem{(b)} Two points on the different incoming segments of the path can be connected by a
differentiable arc which does not intersect the remaining two paths, hence the winding
number on all segments is the same.
\mitem{(c)} Choose points $A',B'$ and $C'$ on the incoming segments of $a,b$ and $c$.
Connect points $A'$ and $C'$ by a differentiable arc $r$ which has exactly one intersection
with path $b$. We can assume that the intersection is transverse and it occurs at point $B'$.
Choose an arc $p$  so that:
\mitem{\Ziffer{1}} $p$ starts at $C'$;
\mitem{\Ziffer{2}} the concatenation of paths $r$ and $p$ satisfies the conditions of Remark
\Ref I.\#.1.9(?);
\mitem{\Ziffer{3}} $p$ intersects paths $a,b$ and $c$ at $\alpha, \beta$ and $\gamma$ many
intersections respectively, all of which are transversal.
\notem
The parity of the winding numbers of points $C', B'$ and $A'$ determining strong
shelteredness is the same as the parity of $\alpha-\beta, \gamma-\alpha+1$ and
$\beta-\gamma +2$, respectively. This follows from Remark \Ref I.\#.1.9(?) using
the path $p$ combined with appropriate segments of $r$ when necessary. Hence the
sum of the winding numbers in question is odd, implying that either one or three of
the winding numbers are odd. Consequently, either one or three of the incoming segments
at $P$ are strongly sheltered.\qed
\2{\bf\Setc{3} Lemma: }
{\rm (Also observed by [Minc], cf.\ his Prop.2.6):}\neuzl
Let $a,b$ and $c$ be any three PL-arcs in the plane from $A$ to $B$ which have
only general position intersections. Then there exists
a strongly sheltered path from $A$ to $B$ within the sum
of traces of paths $a$, $b$ and $c$.
\3{\it  Proof: }
Transverse intersection between two different
paths is the only type of intersection phenomenon that can
occur. Thus the union of traces,
when treated as a graph, has (with the exceptions
of the points $A$ and $B$) only vertices of valency four coming
from the intersections of our paths. Lemma
\Ref I.\#.2.2(?)(a) implies
that from the four incoming edges for each of these vertices exactly two of them are
strongly sheltered and that the sheltered path switches between the two intersecting
paths at the vertex. Furthermore, either one or three of the incoming edges are strongly
sheltered at $A$ and $B$. Thus one can trace a strongly sheltered path starting at $A$.
Such path either ends at  $B$ or comes back to $A$. In the latter case a third edge
from $A$ is also strongly sheltered. A strongly sheltered path starting at such vertex
$A$ can only end at $B$. \qed
\2{\bf\Setc{4} Lemma: }
Let $a,b$ and $c$ be any three PL-paths  in the plane from $A$ to $B$ which have only
general position intersections and self-intersections. Then there exists
a strongly sheltered path from $A$ to $B$ within the sum
of traces of $a$, $b$ \& $c$.
\3{\it Proof: }
In difference to Lemma \Ref I.\#.2.3(?) we now have the additional situation
that a path can intersect itself. Lemma
\Ref I.\#.2.2(?)(b)  implies that either all four or none of incoming edges of such self
intersection are strongly sheltered. Thus the subgraph of strongly sheltered edges is
a graph which has only vertices of even valency with vertices $A$ and $B$ being the
only exceptions by having an odd valency. Similarly as in Lemma \Ref I.\#.2.3(?) we
can find a strongly sheltered path from $A$ to $B$.\qed
\3{\bf\Setc{5}
How a sheltered middle path might solve the three-colored sphere problem: }
Given a map $f : S^2 \longrightarrow \Kr R^2$
where $S^2$ is three-colored as in Question \Ref I.\#.1.2(?), the images of the three
bi-colored meridians are paths that run from $P:= f(N)$ to $Q:= f(S)$. We claim that
the points that are sheltered in the sense of Definition
\Ref I.\#.2.1(?)(b) contain all three colors. Such sheltered points
fall into two categories.

\mitem{(a)} The points which lie on the images of at least two meridians. These points
obviously have all three colors.
\mitem{(b)} The points which lie on the image of one meridian and have a
nonzero winding number with respect to the other two meridians. Choose any such point $x$.
The containment in the image of one meridian implies the presence of the appropriate two colors
at $x$. Considering the third color, note that the image of the correspondingly colored region
on the sphere induces a nullhomotopy of the loop obtained as a concatenation of the other
two meridians. Since $x$ has a nonzero winding number with respect to the other two
meridians any such nullhomotopy contains $x$ in its image (see Definition \Ref I.\#.1.8(?))
hence $x$ also has the third color.
\smallskip\notem
We conclude that all sheltered points have all three colors. Consequently,
wherever the sheltered middle path comes out  as an honest path it is a three-colored
path in the sense of Question \Ref I.\#.1.2(?) and solves the three-colored
sphere-problem in that case.
\initGraphik(1)14.5x15cm
\topinsert\bildbox14.5x12cm{{\bf Figure b} (belonging to
Sect.\Ref I.\#.3.?(?)): \neuzl         
This figure shows that three continuous paths in the plane can
intersect so that the sheltered middle path is
a topologist's sine curve. The construction is naturally a fractal-like
iteration, and only the very few first steps can be pictured.
The construction principle is explained in
\Ref I.\#.3.1(?)--\Ref I.\#.3.10(?)\BoxtextEnd.}\endinsert
%
%
\1{\bf\Paragr{3} . } The general case: an example where the sheltered middle path of  three
non-PL-paths  results in a topologist's sine curve}
\3In this section we present an example of how three paths
in the plane, connecting the same pair of distinct points,
might be placed in such a way that the points that are
sheltered by them in the sense of Definition \Ref I.\#.2.1(?)(b)
yield a topologist's sine-curve but not a genuine path.
The construction is sketched in \Ref I.\#.3.0(Fig.b). It is a plane
drawing which pictures the union of the traces
of our three paths. The construction requires an infinite iteration, hence only the first finitely
many steps can be drawn. Since the definitions near
the regions where our sheltered middle path (the path
in bold on the diagram) turns around (near the left and
right boundary of our diagram) need to be different
from the region in between, the first few steps that we could draw
hardly unveil the construction principle.

Most of this section is devoted to a description of our example. Since the paths are
relatively complicated we shall describe the individual oriented sections of
their traces (\Ref I.\#.3.1(?)--\Ref I.\#.3.10(?)). Afterwards we shall assemble these traces in a continuous
path (\Ref I.\#.3.12(?)--\Ref I.\#.3.14(?)) and eventually compute the appropriate winding
numbers (\Ref I.\#.3.15(?)--\Ref I.\#.3.16(?))
which will confirm that the sheltered path is indeed the bold path of
\Ref I.\#.3.0(Fig.b).
We shall begin by assigning phrases to some of the elements of our construction and the
diagram that will be needed for a complete description.
\3{\bf\Setc{1} Notations used in \Ref I.\#.3.0(Fig.b): }
\bultem
  $P$, $Q$ and $T$ are three points in the plane.
\bultem $[T,Q]$ denotes the straight line segment between these points,
   as it is marked in \Ref I.\#.3.0(Fig.b). It will be also called {\it``the straight line
   segment"} or the {\it``central line"}.
   Within a topologist's sine curve to be defined
   in this section it will serve as the accumulation line.
\bultem $a$, $b$ and $c$ are the names of our three paths. When they show
   up in the figure, they are to be understood in this way,
   that the corresponding segment belongs to our path
   $a$, $b$ or $c$ respectively. Since the paths are highly
   self-intersecting  and intersecting each other along the straight
   line segment, their definition can hardly be understood
   from these labels, but will be explained in
   \Ref I.\#.3.13(?)/\Ref I.\#.3.14(?).
\bultem ``$x$'' and ``$y$'' indicate areas where \Ref I.\#.3.0(Fig.b) is not drawn
correctly for the reasons of scale. Appropriate areas should be changed
   according to the corresponding subfigure of \Ref I.\#.3.3(Fig.c).

\3{\bf\Setc{2} Remark: }
It is easy to verify by counting the winding-number that the bold line from
\Ref I.\#.3.0(Fig.b) is the sheltered path according to
Definition \Ref I.\#.2.1(?)(b).
More explanation will follow in \Ref I.\#.3.15(?).
\par
Furthermore, we shall explain below that
the construction principle can be iterated infinitely many times
 by squeezing in smaller and smaller loops. When
the same intersection patterns as those shown in these first steps are used,
 the line in bold will continue to return forward and backward from
the left boundary-near region of the figure to right one. The
result will be a topologist's sine curve converging to the straight
line segment (cf.\ Remark \Ref I.\#.3.8(?)).
\initGraphik(0)14.5x14.5 cm
\mitBeschriftung
\setbox0=\hbox{\epsfysize=3truecm\epsfbox{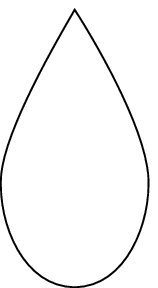}}%
\putbox nach (1.5,7.) dieKoo (m,m) vonBox0%
\setbox1=\hbox{\bleibklein A straight loop}%
\putbox nach (1.5,4.5) dieKoo (m,u) vonBox1%
\setbox0=\hbox{\epsfysize=3truecm\epsfbox{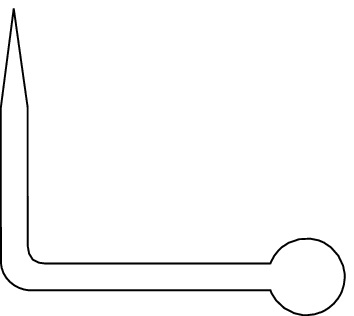}}%
\putbox nach (2.,11.5) dieKoo (m,m) vonBox0%
\setbox1=\hbox{\bleibklein A standard turning loop}%
\putbox nach (2.,9.) dieKoo (m,u) vonBox1%
\setbox0=\hbox{\epsfysize=3truecm\epsfbox{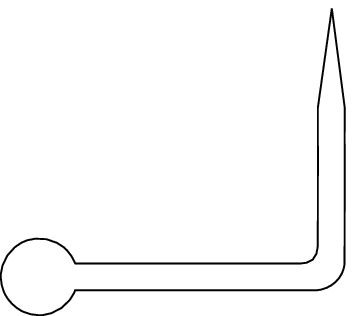}}%
\putbox nach (2,2.5) dieKoo (m,m) vonBox0%
\setbox1=\hbox{\bleibklein A non-standard turning loop}%
\putbox nach (2,0.) dieKoo (m,u) vonBox1%
\setbox0=\hbox{\epsfysize=3truecm\epsfbox{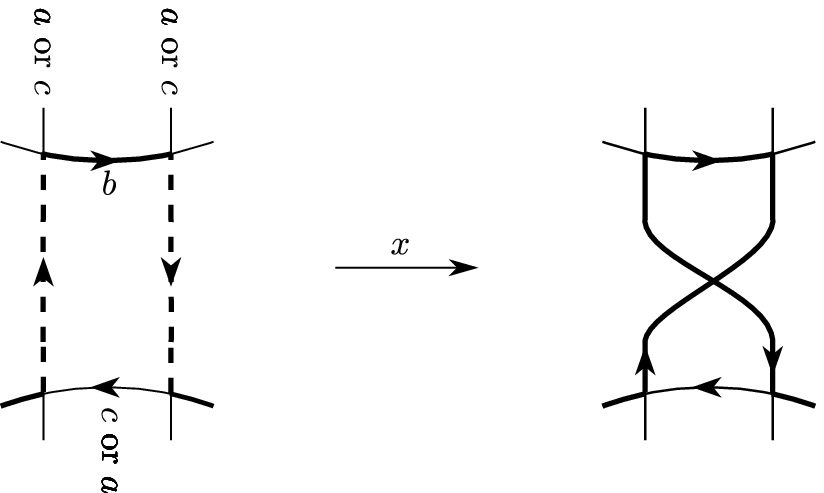}}%
\putbox nach (9,11.5) dieKoo (m,m) vonBox0%
\setbox1=\hbox{\bleibklein ``$x$'' indicates that two vertical strands should actually
cross each other}%
\putbox nach (9,9.) dieKoo (m,u) vonBox1%
\setbox0=\hbox{\epsfysize=3truecm\epsfbox{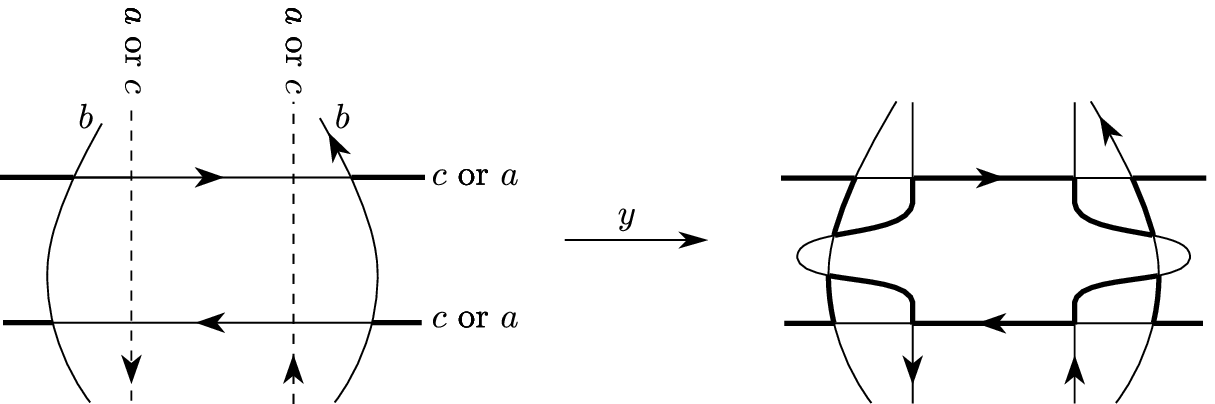}}%
\putbox nach (9.,7.) dieKoo (m,m) vonBox0%
\setbox1=\vbox{\hsize=8truecm\bleibklein\noi``$y$'' indicates that a crossing of a pair of
vertical and a pair of
horizontal strands is to be replaced by a ``clover-crossing''.}%
\putbox nach (9,4.5) dieKoo (m,u) vonBox1%
\setbox0=\hbox{\epsfysize=3truecm\epsfbox{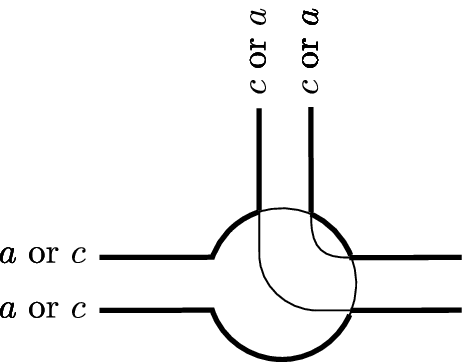}}%
\putbox nach (9,2.5) dieKoo (m,m) vonBox0%
\setbox1=\vbox{\hsize=6truecm\bleibklein\noi Here an arrangement of two standard turning
lops is shown that we call ``hooked''}%
\putbox nach (9,0.) dieKoo (m,u) vonBox1%
\soweitdieBeschriftung
\topinsert\bildbox14.5x14.5cm{{\bf Figure c} (belonging to
\Ref I.\#.3.3(?)--\Ref I.\#.3.10(?)):      
Supporting
     diagrams for \Ref I.\#.3.0(Fig.b)\BoxtextEnd.}\endinsert
\3{\bf\Setc{3} Definition of various kinds of loops: }
Some parts of our paths $a$, $b$ and $c$ will be called ``loops".
According to their geometric shape we shall distinguish between {\it``straight  loops"}\/
and ``turning loops". The turning loops are then again distinguished into
``standard turning" and ``non-standard turning". We decided  to call  the kind which
is more often used in the diagram {\it``standard turning",}\/ and
the other kind {\it``non-standard turning"}. \Ref I.\#.3.3(Fig.c) explains which of the names
is assigned to which kind of loop. The loops
which have the geometric shape of the drop will be called ``straight",
the other ones ``turning". All straight loops will
be labeled ``$b$" and all standard turning loops will be labeled either by ``$a$" or by ``$c$".
The starting and  ending point of each such loop-path is called the {\it``vertex"}\/.
All vertices are on the central line. Considering turning loops  we shall be also talking
about the {\it``vertical part"}\/ and
the {\it``horizontal part"}. Between the vertical and the horizontal part there is
the {\it``turning region",}\/ and opposed to the vertex there is
{\it``the turning-around disk"}.
\3{\bf\Setc{4} The biggest elements of the diagram: }
Looking closer at
\Ref I.\#.3.0(Fig.b) there is no loop turning in the non-standard way therein.
The only part of the diagram which looks approximately like this is
actually not a loop, because it is not a closed curve. It is rather a curve which
connects two nearby different points $P$ and $T$ in a very
indirect way. This line is labeled ``$c$" and will  henceforth be called
{\it``the open loop"}. It is one of the exceptions that do not go into the
fractal-like iteration. The other exception is the leftmost $b$-loop. All
other straight- and standard-turning loops are in principle already a part of
the fractal-like iteration. In particular, the standard-turning loops are  put into
{\it``generations"}. The biggest standard-turning loop in the drawing is labeled by
$a$. It has its vertex at the point $Q$. Similarly as for the other standard-turning
loops of the fractal-like iteration,
most of the vertical part of this loop falls into the inner region of a
straight loop. (Actually, all other loops with vertex at $Q$ have their entire vertical part
within the inner region of the leftmost straight loop.) For this single loop of
first generation the horizontal part
is almost as long as the central line and the turning around disk intersects
with one of the strands of the vertical part of the single open loop.
\3{\bf\Setc{5}  The elements of the second generation of the diagram: }
There are two standard-turning loops of the
second generation in \Ref I.\#.3.0(Fig.b), each approximately half the
size of the single standard-turning loop of the first generation. The size comparison relates
to their horizontal as well as their vertical parts although the horizontal size may become
smaller with respect to the vertical size in later generations. The one on the right hand side
is intersecting with the open loop in a similar way as the
single standard-turning loop of the first generation. This fact will hold for every generation:
the rightmost standard-turning loop of every
generation is always intersecting at its turning-around disk with the left
vertical strand of the open loop. The left standard-turning loop of the
second generation is intersecting with the leftmost $b$-loop. It has its
vertical part inside the leftmost $b$-loop in such a way that inside the
$b$-loop there is no intersection with the vertical part of the single
standard-turning loop of the first generation. This fact will hold for every generation:
the leftmost standard-turning loop will always have its vertical
part in the leftmost $b$-loop, not intersecting the vertical parts of the
turning loops of the prior generations.

All other $b$-loops have another pattern of intersection of their turning loops in the
interior. The $b$-loop at the middle of the central line covering the vertical part of
the right hand turning loop of the second generation is the only
$b$-loop which shows a bit of the general behavior in the course of the
fractal-like iteration. Considering this general behavior: each
vertical part of a turning loop of a new generation is in a similar way
placed in the interior of a straight $b$-loop. ``Similar way" means that the
intersection of the vertical part of the turning loop with the line of the
$b$-loop is above the turning region, but still below the level where
all horizontal parts of the loops of the forthcoming generations intersect this
arrangement of the vertical part of a turning loop.  Each $b$-loop partially contains the
vertical part of exactly one turning loop. It should
also be mentioned, that all those loops of forthcoming generations, which by
their general placement do not intersect with the open loop, are in a
similar way as the left turning loop of the second generation ``hooked" to
the next loop of the same generation on the right. The phrase {\it``hooked"}\/ is
explained in \Ref I.\#.3.3(Fig.c). Recall that our construction is a
plane arrangement, thus ``hooked'' just means that
the turning part of one of the turning loop is
placed in the interior of the turning-around disk of another turning loop,
and that the natural intersections occur in this pattern.
\3{\bf\Setc{6} On the meaning of $y$-regions: }
In the middle of the
$b$-loop that is placed in the center of the segment $[Q,T]$ there is a letter
``$y$" which, similarly as ``$x$",  acknowledges that the given
diagram \Ref I.\#.3.0(Fig.b) is actually incorrectly drawn for reasons of scale.
\Ref I.\#.3.3(Fig.c) explains how the diagram should have been
drawn. The rule is as follows and it applies to all situations where the
horizontal parts of turning loops of higher generations intersect the
vertical parts of a turning loop of a prior generation. It always happens in
the interior of a $b$-loop  where an intersection of two parallel
horizontal and two parallel vertical strands occur. Such intersection situation is to be
replaced by a {\it``clover figure",}\/ as shown in \Ref I.\#.3.3(Fig.c). Note that
with the exception of the leftmost straight loop all other straight loops
have only two vertical strands inside. The letter $y$ places a
clover-figure, which means that wherever the horizontal strands of loops of
higher generations cross the arrangement of a $b$-loop and vertical strands of
an $a$- or $c$-loop inside, the (in general) vertical strands turn between the
two parallel intersecting horizontal strands outside,
intersect the strands of the surrounding $b$-loop, but
then turn around again and come back into the interior of the $b$-loop, and go
again into a vertical position before they intersect the second of the two
parallel horizontal intersecting strands. Although there
is only one letter $y$ in the current figure, the clover figure will
 appear infinitely many times to the two $a$-strands
in the middle $b$-loop on their way going upwards
because there will be infinitely many new generations of loops
between the (currently last drawn) turning loops of the third generation,
and the line $[Q,T]$. Such clover-figures will always be placed inside all
$b$-loops granting them infinitely many clover figures, except inside the leftmost one which,
as we mentioned before, is an exception.

The role of the clover figure is to prevent the sheltered path from reaching the central
line $[Q,T]$ via the vertical strands of the standard turning loops as that would prevent
the sheltered path from being a topologist's sine curve with limit line $[Q,T]$.

\3{\bf\Setc{7}  The elements of the third generation of the diagram:}
We have already discussed how the middle one of the
three turning loops of the third generation intersects the middle
straight loop. Actually, we have discussed all other aspects of the three pictured
turning loops of this generation. To complete the description of the iteration
from the aspect of how the diagram looks we will briefly describe how the straight and
turning loops of the higher generation are put inside this figure.

As already mentioned,
with the exception of the leftmost $b$-loop the vertical parts of turning loops are never
placed inside existing $b$-loops.  They are always located  at new places
between the existing vertical parts of the loops of the prior generation.
Each of them receives its own $b$-loop to cover 80\% of the vertical part.
When putting a new generation of turning loops
into the figure, there may be more than one standard turning loop with its
vertical part placed between the existing vertical parts of loops of prior generation.
Each new generation
is placed in such a way that the horizontal parts together are hooked
and give a through connection between the vertical part which sits
in the leftmost exceptional $b$-loop, and the left line of the vertical
part of the open $c$-loop at the right hand side.

Naturally, with increasing
generation index the loops have to become small. The
vertical parts of the turning loops become small because the horizontal parts
are placed higher with every new generation. The horizontal
parts of the turning loops are becoming smaller by the demand of placing
at least one of them between the vertical
parts of the existing $b$-loops for each new generation. This causes the length
of the horizontal parts of the loops of each new generation
to have essentially  half the length with respect to the loops of the preceding generation.
Word ``essentially" suggests that  there is
an uncertainty of factor $2$ in the length-calculation. This factor of the turning loops
depends on whether the horizontal part has or has not to cross
the vertical part of a loop of the prior generation according to the requirements below.
However, since we make sure that within each new generation
each turning loop crosses at most one loop of a prior generation,
these factors $2$ do not accumulate. There will be at most one
such factor $2$ and it does not prevent the lengths of the horizontal
parts from tending to zero either. Whether we shall place one or two
new turning loops of the new generation between the existing turning loops and whether
these loops will be labeled by $a$ or by $c$, has to be arranged so that the
following rules are kept:
\bultem the loops of a new generation which intersect the vertical part
   of a turning loop and the surrounding $b$-loop of a prior
   generation have to have a different
   label than the strands of the loop of the old generation that are
   inside the $b$-loops in the vertical position, i.e., inside a $b$-loop
   only $c$ and $a$ are allowed to intersect but not $a$ with
    $a$ or $c$ with $c$.
\bultem The rightmost turning loop of each new generation (the one intersecting
   the open $c$-loop) has to be ``of type $a$", i.e., has to have the
   label ``$a$".
\bultem Every two loops which are hooked have to be of different type.
\3{\bf\Setc{8}  Remark \Ref I.\#.3.2(?) continued: }
Based on these demands it should be  clear that the construction
can be continued to infinity maintaining these demands and adding
smaller and smaller loops to our diagram. Observe that each
new generation of turning loops requires a new generation of
$b$-loops, each covering approximately 80\% of the vertical part
of a turning loop of the new generation. The only
exception is the leftmost loop of each new generation, which goes
always into the leftmost old $b$-loop of the first generation.
The situation for the first three generations is drawn
in the \Ref I.\#.3.0(Fig.b).
\3{\bf\Setc{9} Definition of the orientation of loops: }
We now clarify how the orientation of loops as suggested in \Ref I.\#.3.0(Fig.b) coincides
with the orientation of the segment of the corresponding path
$a$, $b$ or $c$ that traverses this loop as a part of its trace.
The straight line should not be regarded as oriented
since it is in both directions run through by our paths. For this reason the calculation of
the winding numbers will avoid it.
However,  all other elements of our figure are in the trace
of one (or at most two at the intersections) of our paths. Hence they should be regarded as
oriented which is indicated by our \Ref I.\#.3.0(Fig.b). To repeat
the rules and to clarify the conventions for the
parts of the iteration that could not have been drawn, we
summarize them within the following list:
\bultem The three lines that are emanating from the point $P$ are oriented away
   from the point $P$. This automatically implies that the open $c$-loop
   is clockwise oriented.
\bultem All straight loops (i.e., all $b$-loops) are anticlockwise oriented.
\bultem All turning loops (i.e.,  all $a$-loops and all $c$-loops)
   are in principle  clockwise oriented. However,  the presence of $x$-regions effects such
   orientation which is explained in
   \Ref I.\#.3.10(?).
\3{\bf\Setc{10} On the meaning of $x$-regions: }
``$x$" (similarly to ``$y$") means that the drawing of \Ref I.\#.3.0(Fig.b) is not correctly
drawn but
needs to be changed. Such letter is to be added to every (in its
generation not leftmost) standard turning loop in the lower part of the vertical section
where the loop is already above the turning-around disk to which
it is hooked, and still below the intersection of the straight $b$-loop which
is supposed to cover 80\% of its vertical part. Letter
$x$ in principle means that in this part the two vertical strands of our
loop are not (as they are pictured in the current drawing) parallel but
actually cross each other. A consequence of this crossing is that the
orientation of the loop can be only clockwise in the horizontal part, the
turning region, and the turning-around disk, but in the majority of the vertical
part (due to this crossing) the orientation will naturally be
anticlockwise. Therefore the orientation is not contradictory
as it might appear at first sight when looking at \Ref I.\#.3.0(Fig.b).
\3{\bf\Setc{11} A note on continuity : }\neuzl
Recall (cf.\ the introductory paragraph of  this section) 
that
\Ref I.\#.3.0(Fig.b) together with the explanations given
in the preceding paragraphs on how to iterate its construction
was set up to be part of a definition of three paths $a$, $b$ and $c$ in the plane
with the following properties:
\bultem each of the three paths connects points $P$ and $Q$;
\bultem the three paths induce the bold path of \Ref I.\#.3.0(Fig.b) as a sheltered middle
path (according to Definition \Ref I.\#.2.1(?)(b)) which is a topologist's sine curve.
\notem
The description of \Ref I.\#.3.0(Fig.b) and the
principle according to which the iteration of the construction
should be continued is by now completed. We proceed with an explanation on
how the continuous
paths $a$, $b$ and $c$ can be defined so that they traverse
infinitely many corresponding
loops. A concatenation of countably many paths of appropriate diameters
can be performed in various ways. The following lemma is adapted to our situation.

\2{\bf\Setc{12} Lemma on infinite concatenation of paths: }
Suppose $\alpha\colon [0,1]\to \Kr R^2$ is an arc between points $A$ and $B$,
$\{x_i\}_{i\in \Kr N}$ is a dense subset of $\alpha((0,1))$ and
$\{\alpha_i\colon[0,1]\to \Kr R^2\}_{i\in \Kr N}$ is a collection of loops based at
$x_i$ such that the sets $\alpha([0,1]),\alpha_1((0,1)), \alpha_2((0,1)), \ldots$ are
pairwise disjoint and such that $\lim_{i\to\infty}{diam}(\alpha_i([0,1]))=0$. Then
there exists a path $\beta\colon [0,1]\to \Kr R^2$ from $A$ to $B$
with $$\beta([0,1])\subset \alpha([0,1])\cup \bigcup_{i\in \Kr N}\alpha_i((0,1)),$$
such that $\beta$ traverses each loop $\alpha_i$ exactly once (preserving the
orientation of $\alpha_i$), i.e.,

  \bultem for each $i\in \Kr N$ there exists $(a_i,b_i)\subset (0,1)$
  so that $\beta(a_i+(b_i-a_i)\mal t)=\alpha_i(t), \forall t\in(0,1)$;
  \bultem for each $i\in \Kr N$ and every $x\in \alpha_i((0,1))$ the
  preimage $\beta^{-1}(\{x\})$ is exactly one point.

\3{\it Proof: } Path $\beta $ will be obtained as a limit path of a uniformly
convergent sequence of paths $\beta_n \colon [0,1]\to \Kr R ^2$. The iterative
construction of paths $\beta_n$ will imitate the construction of the standard ternary Cantor
set $C$, which is obtained by inductively removing the middle thirds of the
intervals. Accordingly we denote by $C_n$ an appropriate inductive step of the
construction of $C$ so that $\cap _{n=0}^\infty C_n=C$, i.e., $C_0=[0,1]$,
$C_1=[0,1/3]\cup[2/3,1]$,{\tt~...}

Path $\beta_1$ is obtained from path $\alpha$ by inserting loop $\alpha_1$ in the
appropriate spot. We divide the interval $[0,1]$ into three equal parts and define:
$$
\beta_1(t) := \left\{\vcenter{\halign{\vrule width0pt height 8pt depth 8pt$#$\hfil&~~for~$#$\hfil\cr
\alpha(3\mal t\mal x_i) & t\in [0,{1\over3}],\cr
\alpha_1(3\cdot(t-{1\over3})) & t\in [{1\over3},{2\over3}],\cr
\alpha(3\mal(t-{2\over3})\cdot(1-x_i)+x_i)&t\in [{2\over3},1].\cr}}
\right.$$
\notem
In other words, $\beta_1|_{[0,1/3]}$ is the reparameterized path $\alpha|_{[0,x_1]}$,
$\beta_1|_{[2/3,1]}$ is the reparameterized path $\alpha|_{[x_1,1]}$ and
$\beta_1|_{[1/3,2/3]}$ is the reparameterized loop $\alpha_1$. Consequently
$(a_1,b_1)=(1/3,2/3)$ and $\beta_1|_{C_1}$ essentially represents loop $\alpha$,
i.e. the concatenation of $\beta_1|_{[0,1/3]}$ and $\beta_1|_{[1/3,2/3]}$ is the reparameterized
path $\alpha$.

In a similar way we can define $\beta_2$ as a path obtained from $\beta_1$ by
inserting loops $\alpha_2$ and $\alpha_3$ in the appropriate spot so that
$(a_2,b_2)=(1/9,2/9)$ and $(a_3,b_3)=(7/9,8/9)$. Note that the properties of paths
$\alpha_i$ are preserved under every permutation of indices $i$, hence we can assume
that $0<x_2<x_1<x_3<1$. Furthermore, $\beta_2|_{C_2}$ essentially represents the loop
$\alpha$ in a similar way as $\beta_1|_{C_1}$ essentially represents loop $\alpha$.

We proceed by induction. By redefining $\beta_n$ on $C_{n}$ so that the middle thirds
of the intervals are mapped by ever smaller loops $\alpha_i$ and the remainder (which
essentially represents path $\alpha$) is reparameterized, we obtain $\beta_{n+1}.$
Note that the sequence $(\beta_i)_i$ is uniformly convergent as diameters of
$\alpha_i([0,1])$ as well as diameters of components of $C_i$ tend to $0$. Hence we
obtain a continuous path $\beta$ satisfying the required conditions. \qed

\3{\bf\Setc{13} Interpretation of Lemma \Ref I.\#.3.12(?): }

Consider the notation of Lemma \Ref I.\#.3.12(?). Note that the map
$\beta|_{[0,1]-C}$ essentially consists of loops $\alpha_i$ and $\beta|_C$ represents
the path $\alpha$, i.e., there exists a surjection $g\colon C \to [0,1]$ so that
$\beta=\alpha \circ g$.
The behavior of a path $\beta$ when constructed in this way is described by the following
statement:
     {\it``Path $\beta$ moves on a Cantor set over a line segment $\alpha([0,1])$ and performs
     excursions $\alpha_i([0,1])$ on all intervals not belonging to the Cantor set"}.

\3{\bf\Setc{14} Definition of paths $a$, $b$ and $c$: }\neuzl
The definition of the paths $a$, $b$ and $c$ is based on the construction principle
explained in \Ref I.\#.3.12(?). All three of them start at $P$ and end at $Q$.

Definition of the path $a$:
\bultem starts at $P$;
\bultem runs over the short curve to $T$;
\bultem moves over the straight line segment on a Cantor set to $Q$ while taking
   excursions on all the correspondingly marked $a$-loops not belonging to the
   Cantor-set.\notem

Definition of the path $b$:
\bultem starts at $P$;
\bultem runs over the accordingly marked long curve directly to the point $Q$;
\bultem moves on a Cantor set over the central line to the point $T$ while
   taking excursions over all $b$-loops on the intervals that do not belong
   to this Cantor set;
\bultem runs straight without further excursions  back from $T$ to $Q$.\notem

Definition of the path $c$:
\bultem  starts at $P$;
\bultem  runs over the single open non-standard turning loop to $T$;
\bultem moves on a Cantor set over the central line
   to the point $Q$ taking excursions over all $c$-loops on the
   intervals that do not belong to the Cantor-set.\notem

Thus we have defined continuous paths
containing all the correspondingly marked loops
and other segments of \Ref I.\#.3.0(Fig.b) and all infinitely many loops that were to be
added according to the iteration described in
\Ref I.\#.3.4(?)--\Ref I.\#.3.10(?).
\3{\bf\Setc{15} Comments on the resulting sheltered path: }\neuzl
We already mentioned that the sheltered path is highlighted
by boldness in \Ref I.\#.3.0(Fig.b). At
the critical, most details requiring phases, the information about boldness
should be taken from the corresponding subfigures of \Ref I.\#.3.3(Fig.c).
This in particular affects the hooked regions where everything
is becoming very small and  the regions
where the corrections via letters $x$ and
$y$ are implemented. The shelteredness in all above mentioned regions
behaves analogously and can be computed by a procedure described in Remark
\Ref I.\#.1.9(?).

With the exception of the central line all parts of the diagram
are tame. Outside the central line each of
the drawn loops is only once traversed by
one segment of the path whose label it carries. The visible
intersection points are the only double points. Furthermore,
there will be no additional insertions below the three generations
of loops that we have drawn.

In order to verify that the sheltered path has been correctly marked we use the criterion of
Remark \Ref I.\#.1.9(?). Given any point on a bolded section choose an arc according to
Remark \Ref I.\#.1.9(?). In fact we can choose a straight line segment transversely
intersecting
the pictured path-system and avoiding the central line.
Given such arc the intersection numbers and consequently the winding numbers can be obtained.
Note that we have defined the appropriate  winding numbers up to a factor $-1$.
\3{\bf\Setc{16} Conclusion for the resulting sheltered path:}
The following list describes the resulting sheltered path which starts at $P$.
\bultem The path runs from $P$ through the open $c$-loop until
its intersection with the turning-around disk of the single standard-turning
loop of the first generation. This section has winding number $1$ with respect to the
other two paths by Remark \Ref I.\#.1.9(?).
\bultem The second part runs forward and backward over the horizontal
parts of the single standard-turning loop of the first generation labeled by $a$,
turning around at the straight loop labeled $b$. The appropriate winding number of
this section is $1$. In particular, the winding number with respect to $b$
and $c$ is $1$  on the contained part of the standard turning loop labeled $a$.
Similarly,  the winding number with respect to $a$ and $c$ is $1$  on the contained
part of the straight loop labeled $b$. All other paths approaching the intersections
in question have winding number $0$.
\bultem The path switches between the rightmost standard turning loops  of the first
and the second generation along the path labeled $c$, connecting their turning-around regions.
The corresponding winding number is $1$.
\bultem The path continues along the two hooked horizontal parts of the
loops of the second generation and reaches the turning-around region of the rightmost
standard turning loop. In particular, it runs forward and backward along the horizontal
parts, turning around at the straight loop labeled $b$, switching between the loops of
the second generation and making a small excursion in the $x$ region according to
Figure $c$. The corresponding winding numbers are $1$ whose calculation should also
include the information of \Ref I.\#.3.3(Fig.c).
\bultem The path continues to switch to the rightmost standard turning loop of the
next generation and to run forward and backwards on the hooked horizontal
parts of the turning loops of every generation.  Small diversions occur
in the hooking-regions and in the clover-crossings wherever the vertical parts of our
turning loops have to be crossed.

Since the running backward
and forward always occurs between the left vertical segment of the open
$c$-loop on the right hand side of the picture and the right strand of the
leftmost $b$-loop on the left hand side of the picture, the obtained path is a topologist's
sine curve accumulating to the straight line.

\1{\bf Acknowledgment}

\3 This research was supported by the Polish--Slovenian  grants BI-PL 2008-2009-010
and 2010-2011-001,
the first and the third author by
the ARRS program P1-0292-0101 and
project J1-2057-0101
and the second and the fourth author by the MNiSW grant N200100831/0524.
We thank the referee for comments and suggestions.

\bigskip\4%
\mitem{[Bgt]} BOGATYI, S.A.:~~``The topological Helly theorem" (in Russian),
Fund.\ Prikl.\ Mat., Vol.\ 8 (2002), pp.\ 365--405;
\mitem{[Breen]} BREEN, Marilyn:~~ ``A Helly-type theorem
for intersections of compact connected sets in the plane",  Geom.\ Dedicata,
Vol.\  71  (1998),  no.\ 2, pp.\ 111--117;
\mitem{[DGK]} DANZER, L., GR\"UNBAUM, B.\ and KLEE, V.:~~``Helly's theorem and its
   relatives. {\bf In} ``Convexity", edited by
   Victor L.\ Klee, Proc.\
   Sympos.\ Pure Math. Vol.\ VII, Amer.\ Math.\ Soc.,
   Providence (R.I.), U.S.A., 1963, pp.\ 101--180;
\mitem{[Eck]}  ECKHOFF, J\"urgen:~~``Helly, Radon, and Carath\'eodory type theorems".
{\bf In:} ``Handbook of convex geometry (Vol.\ A)", edited by
    P.\ M.\ Gruber and J.\ M.\ Wills,
    North-Holland Publishing Co., Amsterdam, 1993, pp.\ 389--448;
\mitem{[Farb]} FARB, Benson:~~``%
Group actions and Helly's theorem",
Adv.\ Math., Vol.\ 222 (2009), no.\ 5, pp.\ 1574--1588;
\mitem{[Flstd]} FL{\O}YSTAD, Gunnar:~~``The colorful Helly theorem and
colorful resolutions of ideals", in press at J.\ Pure Appl.\ Algebra,
Vol.\ 215 (2011), no.\ 6, pp.\ 1255--1262;
\mitem{[Hlly1]} HELLY, Eduard:~~``\"Uber Mengen konvexer K\"orper mit
    gemeinschaftlichen Punkten",
    Jahresberichte der Deutschen Math.-Ver., Bd.~32 (1923), pp.\ 175--176;
\mitem{[Hlly2]} --- , ``\"Uber Systeme von abgeschlossenen Mengen mit gemeinschaftlichen
Punkten", Mo\-natsh.\ Math.\ Phys., Bd.~37 (1930), pp.\ 281--302;
\mitem{[KaRe]} KARIMOV, U.\ H.\ and REPOV\^S, D.:~~``On topological Helly theorem'',
   Topol.\ Appl., Vol. 153, No.\ 10 (2006), pp.\ 1614-1621;  
\mitem{[KR\^Z]} --- and \^ZELJKO, M.;~~``On the unions and intersections of
   simple connected planar sets",    Monatsh.\ Math.\  145  (2005),
   no.\ 3, pp.\ 239--245;
\mitem{[Minc]} MINC, Piotr:~~``Choosing a sheltered middle path'',
   Topology Appl., Vol.\  153  (2006),  no.\ 10, pp.\ 1622--1629;
\mitem{[Mol1]} MOLN\'AR, J\'ozsef.:~~``\"Uber eine Verallgemeinerung auf die Kugelfl\"ache
   eines topologischen Satzes von Helly'', Acta Math.\ Acad.\ Sci.\ Hugar.,
   Vol.\ 7 (1956), pp.\ 107--108;
\mitem{[Mol2]}
--- ,  ``\"Uber den zweidimensionalen topologischen Satz von Helly"
(in Hungarian),  Mat.\ Lapok, Vol.\  8  (1957), 108--114;
\mitem{[TyVa]} TYMCHATYN, E.D.\ and VALOV, V.:~~``On
intersection of simply connected sets in the plane",
Glas.\ Mat.\ Ser.\ III, Vol.\ 41 (2006), no.\ 1, pp.\ 159--163;
\mitem{[Wng]} WENGER,  Rephael:~~``Helly-type theorems and geometric transversals". {\bf In:} ``Handbook of Discrete and Computational Geometry",
CRC Press Ser.\ Discrete Math.\ Appl.,
edited by Jacob E. Goodman and Joseph O'Rourke,  CRC Press, Boca Raton
(Florida, USA), 1997, pp.\ 63--82.
\medskip
\3\tt\halign{# \hfil&\qquad\qquad#\hfil\cr
Du\^san Repov\^s&Witold Rosicki\cr
Faculty of Education&Institute of Mathematics\cr
\quad and Faculty of Mathematics and Physics&University of Gdansk\cr
University of Ljubljana\cr
P.~O.~Box 2964 &ul.\ Wita Stwosza 57 \cr
Ljubljana 1001&80-952 Gda\'nsk\cr
SLOVENIA&POLAND\cr
Email:&Email: \cr
dusan.repovs@guest.arnes.si&wrosicki@mat.ug.edu.pl\cr
\hbox{~~~}\cr\hbox{~~~}\cr
\^Ziga Virk  &Andreas Zastrow\cr
Faculty of Mathematics&Institute of Mathematics\cr
\quad and Physics&University of Gdansk\cr
University of Ljubljana\cr
Jadranska 19 &ul.\ Wita Stwosza 57 \cr
Ljubljana 1000&80-952 Gda\'nsk     \cr
SLOVENIA&POLAND\cr
Email:&Email: \cr
zigavirk@gmail.com&zastrow@mat.ug.edu.pl\cr}
\ENDE